\documentclass[%
 aip,
 amsmath,amssymb,
 reprint,%
]{revtex4-2}

\usepackage{graphicx}
\usepackage{dcolumn}
\usepackage{bm}

\usepackage[utf8]{inputenc}
\usepackage[T1]{fontenc}
\usepackage{mathptmx}
\usepackage{etoolbox,makecell,hyperref}
\usepackage{amsmath,amsfonts,amssymb,amsthm,amstext,enumitem,comment}
\usepackage[usenames,dvipsnames]{xcolor}

\makeatletter
\def\@email#1#2{%
 \endgroup
 \patchcmd{\titleblock@produce}
  {\frontmatter@RRAPformat}
  {\frontmatter@RRAPformat{\produce@RRAP{*#1\href{mailto:#2}{#2}}}\frontmatter@RRAPformat}
  {}{}
}%
\makeatother
\begin{document}

\newcommand{\R}{\mathbb R}
\newcommand{\N}{\mathbb N}
\newcommand{\T}{\mathbb T}
\newcommand{\ep}{\varepsilon}
\newcommand{\lb}{\lambda}
\newcommand{\W}{\Omega}
\newcommand{\w}{\omega}
\newcommand{\wt}{\omega{\cdot}t}
\newcommand{\ws}{\omega{\cdot}s}
\newcommand{\mA}{\mathcal{A}}
\newcommand{\mB}{\mathcal{B}}
\newcommand{\mC}{\mathcal{C}}
\newcommand{\mD}{\mathcal{D}}
\newcommand{\mI}{\mathcal{I}}
\newcommand{\mJ}{\mathcal{J}}
\newcommand{\mL}{\mathcal{L}}
\newcommand{\mM}{\mathcal{M}}
\newcommand{\mP}{\mathcal{P}}
\newcommand{\mR}{\mathcal{R}}
\newcommand{\mb}{\mathfrak{b}}
\newcommand{\mc}{\mathfrak{c}}
\newcommand{\md}{\mathfrak{d}}
\newcommand{\mf}{\mathfrak{f}}
\newcommand{\mg}{\mathfrak{g}}
\newcommand{\mh}{\mathfrak{h}}
\newcommand{\ml}{\mathfrak{l}}
\newcommand{\muk}{\mathfrak{u}}

\newtheorem{teor}{Theorem}[section]
\newtheorem{lema}[teor]{Lemma}
\newtheorem{prop}[teor]{Proposition}
\newtheorem{coro}[teor]{Corollary}
\theoremstyle{definition}
\newtheorem{defi}[teor]{Definition}
\newtheorem{remark}[teor]{Remark}
\preprint{AIP/123-QED}

\title[Nonautonomous uniform stability or bistability in optical fluorescence]{Nonautonomous scalar concave-convex differential equations: conditions for uniform stability or bistability in a model of optical fluorescence}
\author{J. Due\~{n}as}
\affiliation{Departamento de Matem\'{a}tica Aplicada, Universidad de Va\-lladolid\\ Escuela de Ingenier\'{\i}a Inform\'{a}tica de Valladolid, Pº de Bel\'{e}n, 15, 47011 Valladolid, Spain.}
\affiliation{Instituto de Investigaci\'{o}n en Matem\'{a}ticas (IMUVa), Universidad de Valladolid.}

 \author{C. N\'{u}\~{n}ez}%
\affiliation{Departamento de Matem\'{a}tica Aplicada, Universidad de Va\-lladolid\\ Escuela de Ingenier\'{\i}as Industriales, Prado de la Magdalena 3-5, 47011 Valladolid, Spain.}
\affiliation{Instituto de Investigaci\'{o}n en Matem\'{a}ticas (IMUVa), Universidad de Valladolid.}

\author{R. Obaya}
\affiliation{Departamento de Matem\'{a}tica Aplicada, Universidad de Va\-lladolid\\ Escuela de Ingenier\'{\i}as Industriales, Prado de la Magdalena 3-5, 47011 Valladolid, Spain.}
\affiliation{Instituto de Investigaci\'{o}n en Matem\'{a}ticas (IMUVa), Universidad de Valladolid.}

\email{rafael.obaya@uva.es}
\email{carmen.nunez@uva.es}
\email{jesus.duenas@uva.es}

\date{\today}

\begin{abstract}
The long-term dynamics of a Bonifacio-Lugiato model of optical superfluorescence is investigated. The scalar ordinary differential equation modelling the phenomenon is given by a concave-convex autonomous function of the state variable that is excited by a time-dependent input, $I(t)$. The system's response is described in terms of the dynamical characteristics of the input function, with particular focus on uniform stability or bistability cases. Building on previous published results, the open interval defined by the constant input values for which the equation exhibits uniform stability or bistability is considered, and it is proved that bistability occurs when $I(t)$ lies within this interval. This condition is sufficient but not necessary. Applying nonautonomous bifurcation methods and imposing more restrictive conditions on the variation of $I(t)$ makes it possible to determine the necessary and sufficient conditions for bistability and to prove that the general response is uniform stability when these conditions are not satisfied. Finally, the case of a periodic input that varies on a slow timescale is analyzed using fast-slow system methods to rigorously establish either a uniformly stable or a bistable response.
\end{abstract}

\maketitle

\begin{quotation}
The methods of nonautonomous dynamical systems make it possible the extension of classical results on uniform stability or bistability for the Bonifacio-Lugiato optical superfluorescence model to the case of a time-dependent input. The necessary and sufficient conditions for the presence of bistability are described in terms of the relationship between the variations of the external input and the intrinsic map that describes the scalar model, and it is shown that the absence of bistability corresponds, in general, to the presence of uniform stability.
\end{quotation}
\section{Introduction}
A dynamical system describes the evolution over time of a given phenomenon according to certain laws or rules that govern its behavior. When the laws themselves are explicitly time-dependent, the system is said to be nonautonomous. The concept of nonautonomous dynamical system encompasses a solid mathematical theory with its own methods and tools, typically different from those employed in the study of autonomous models.

In this paper, we investigate the dynamics induced by the well-known Bonifacio-Lugiato model of optical superfluorescence,\cite{bonifaciolugiato1,bonifaciolugiato2,bonifaciolugiato3} excited by a non-negative input which, in our case, is assumed to be time-dependent; namely,
\begin{equation}\label{eq:intro}
 x'=I(t)+g(x)
\end{equation}
with $g(x):=-x-2\,c\,x/(1+x^2)$ for a fixed constant $c>0$. In the line of the results obtained in Refs.~\onlinecite{bonifaciolugiato1,bonifaciolugiato2,bonifaciolugiato3}, we provide precise conditions on $I(t)$ that guarantee a uniformly stable response of the system---exactly a bounded solution, which is hyperbolic attractive---as well as others that imply the bistability of the model---the coexistence of two stable states of this type. In particular, these results explain the transition from uniform stability to uniform bistability in the nonautonomous formulation of the problem. It is known that bistability or, more generally, multistability can be advantageous in certain applications, while posing a drawback in others. In scenarios where it is desirable to maintain strict control over the state of a system subject to noisy perturbations, bistability may present a challenge. Conversely, in other contexts, bistability can offer significant flexibility by enabling transitions between states when the system is subjected to appropriate control.\cite{feudel2014,feudel2018} In the case of optical devices, the occurrence of multiple stable states with different light intensities, along with the ability to transition between them, has applications in optical communications and logic gates.\cite{wieczorek1999,wieczorek2005,ashwin2007}

Since $I(t)\ge 0$ and $g(0)=0$, the region $x\ge 0$ is positively invariant. In fact, the model only makes sense for $x\ge 0$: $x$ is proportional to the light intensity. In addition, the $C^1$ map $g\colon[0,\infty)\to(-\infty,0]$ is strictly convex on $[0,a)$ and strictly concave on $(a,\infty)$ for a certain $a>0$, its derivative $g'$ is strictly concave---a property referred to as d-concavity---on an interval $(a-b,a+b)$ with $b\in(0,a)$, and it satisfies $\lim_{x\to\infty} g(x)=-\infty$.  In the study of a nonautonomous scalar ordinary differential equation $x'=f(t,x)$, the global concavity or convexity of the section maps $x\mapsto f(t,x)$, as well as the global concavity of the derivatives $x\mapsto f_x(t,x)$, limit the possibilities of the global behavior and, consequently, also of the global---nonautonomous---bifurcation diagram for $x'=f(t,x)+\lb$: see, e.g., Refs. \onlinecite{alonsobaya}, \onlinecite{nuobsanz}, \onlinecite{dno1}. In fact, the list of models that respond to concave or d-concave differential equations is large. For instance, in Ref.~\onlinecite{dlo}, one can find an account of scalar ecological models of these types. But, of course, most of the models do not fit this global structure: in a general situation, $f$ is convex with respect to $x$ in some areas of its domain, and convex in the complement areas. This is the case of \eqref{eq:intro}.

The recent Ref.~\onlinecite{dno6} presents a first analysis of nonautonomous scalar ordinary differential equations $x'=f(t,x)$ given by maps $f$ that alternate from convexity to concavity: it considers the case of existence of a smooth curve $t\mapsto a(t)$ such that $x\mapsto f(t,x)$ is concave when $x\geq a(t)$ and convex when $x\leq a(t)$, showing the occurrence of dynamical scenarios that had already been described in Refs.~\onlinecite{alonsobaya,nuobsanz,dno1}, as well as new ones arising from the interplay between concave and convex dynamics, which do not occur when $f$ is purely concave or d-concave. In some cases, the existence of a band around the graph of $a$ on which $f$ is d-concave (i.e., $f_x$ is concave) in $x$ allows to delve deeper into some points of the description of the global dynamics. Applying the results of Ref.~\onlinecite{dno6} to the Bonifacio-Lugiato model \eqref{eq:intro} is  the fundamental tool of this paper: the constant map $t\mapsto a$ plays the role of $a(t)$, and $\R\times[a-b,\,a+b]$ is the d-concavity band.

After the preliminary Section \ref{sec:bonifacio-lugiato}, on which the main physical characteristics of the model and the dynamical formulation required to undertake its study are described, we check that bistability is only possible for \eqref{eq:intro} if it is possible for some constant inputs $\lb$ in the autonomous version $x'=\lb+g(x)$. In turn, this requires $g$ to have only a local minimum $x_1$ on $(0,\,\infty)$, lying on $(0,\,a)$ and taking a value $-\lb_2$, and only a local maximum $x_2$, lying on $(a,\,\infty)$ and taking a value $-\lb_1>-\lb_2$.  This is the case if and only if $c>4$, what we assume from now on. So, $\lb_1$ and $\lb_2$ are the two bifurcation values of $x'=\lb+g(x)$, both of saddle-node type. We prove  bistability when $\lb_1<\inf_{r\in\R}I(r)$ and $\lb_2>\sup_{r\in\R}I(r)$ (i.e., when $0\in \mI_1:=(\lb_1-\inf_{r\in\R}I(r),\,\lb_2-\sup_{r\in\R}I(r))$), and uniform stability if either $\sup_{r\in\R}I(r)<\lb_1$ or $\inf_{r\in\R}I(r)>\lb_2$. Clearly, these cases are far away to exhaust the possibilities unless $I(t)$ is constant. We also establish conditions on $I(t)$ extending these first results to the cases of non-strict inequalities. This is done in Section \ref{sec:bistability_autonomous}.

From this point, we will establish more restrictive conditions on the relation between $I$ and $g$ that guarantee bistability when 0 belongs to an interval larger than $\mI_1$. More precisely, we associate two auxiliary equations to \eqref{eq:intro}: $x'=I(t)+g_-(x)$ of globally (non-strict) concave type, and $x'=I(t)+g_+(x)$ of globally (non-strict) convex type. The map $g_-$ (resp.~$g_+$) is defined as the $C^1$ linear continuation of $g$ outside $[a,\,\infty)$ (resp.~outside $[0,\,a]$). Ref.~\onlinecite{dno6}~establishes the existence of a unique---saddle-node type---nonautonomous bifurcation value $\lb_\mp$ for the parametric problems $x'=\lb+I(t)+g_\mp(x)$: in the concave-linear case, there are two separate hyperbolic solutions for $\lb>\lb_-$ and no bounded solutions if $\lb<\lb_-$, and the situation is symmetrical in the linear-convex case. We check that $\mI_1\subseteq\mI_2:=(\lb_-,\,\lb_+)$, and  establish conditions relating the variation of $I(t)$ on its domain with the variation of $g$ on $[x_1,\,a]$ and on $[a,\,x_2]$ which ensure that bistability occurs if and only if $0\in\mI_2$, while uniform stability is the response if $0\notin \overline \mI_2$. This is the main result of Section \ref{sec:bistability_concaveconvex}, which is completed in Section \ref{sec:precise_estimates} with the analysis of the scope of the conditions required on the variation of $I(t)$ for some specific types of inputs.

The study of bistability continues in Section \ref{sec:bistability_dconcavity}, where some additional conditions on $I(t)$---including that it takes values on a compact subset of the set $(-g(a-b),-g(a+b))$, determined by the d-concavity band, as well as possibly smaller variation of the input map---guarantees that all the hyperbolic solutions lie within the d-concavity band and provides an interval $\mI_3\supseteq\mI_1$ such that bistability occurs for $0\in\mI_3$. Finally, in Section \ref{sec:relxationoscillations}, we follow the theory of Tikhonov and Fenichel\cite{tikhonov1952,fenichel1979} to  investigate the case of a periodic input varying in slow time. In some cases, we find three hyperbolic solutions implying bistability while in others we get the so-called relaxation oscillation, which is the unique (exponentially stable) bounded solution of the equation.

\section{The Bonifacio-Lugiato model}
\label{sec:bonifacio-lugiato}
In this paper, we deal with a model of optical superfluorescence, described in Refs. \onlinecite{bonifaciolugiato1} and \onlinecite{bonifaciolugiato2},
\begin{equation}\label{eq:autonomousbonifacio}
x'=\lb-x-2\,c\,\frac{x}{1+x^2}\,.
\end{equation}
It is a mean-field model of a homogeneously broadened ensemble of two-level atoms driven by a coherent resonant field.

\begin{figure}
     \centering
         \includegraphics[width=\columnwidth]{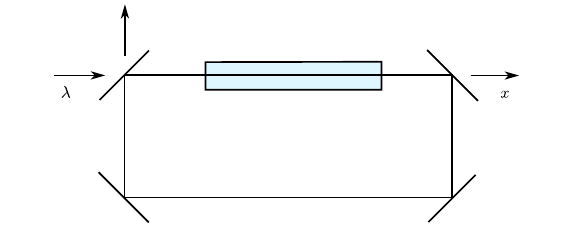}
         \caption{Sketch of a ring cavity formed by four mirrors, containing a medium composed of two-level atoms (in blue).
         The amplitude of the incident field is proportional to the parameter $\lb$, while the amplitude of the transmitted field is proportional to the variable $x$.
         A third arrow represents the reflected field, which is not included in the equation under study.
         The sketch is based on Fig. 1 from Ref.~\onlinecite{bonifaciolugiato3}.}
        \label{fig:bonifaciolugiatointerferometer}
\end{figure}

The system can be conceptualized as a ring cavity containing an ensemble of two-level atoms (see Fig.~\ref{fig:bonifaciolugiatointerferometer}).
It emits a transmitted field whose intensity is proportional to the variable $x$, whose temporal evolution is governed by the equation.
This occurs when the cavity is excited by an incident field with an amplitude proportional to the parameter $\lb$.
Thus $x,\lb\geq0$.
The equilibrium branches of the model are described in Ref.~\onlinecite{bonifaciolugiato1},
and the out-of-equilibrium dynamics of the model is detailed in Section V of Ref.~\onlinecite{bonifaciolugiato2}.
The constant $c>0$ depends on the material used, as it is proportional to the ratio of the incoherent transverse atomic relaxation rate to the cooperative damping rate of pure superfluorescence.

A natural question to consider in the study of this model is how the dynamics changes as $\lb$ varies.
This leads to a bifurcation problem in $\lb$.
In their papers,\cite{bonifaciolugiato1,bonifaciolugiato2} Bonifacio and Lugiato demonstrate that the condition $c>4$ is both necessary and sufficient for the existence of an interval of positive values of $\lb$ for which the autonomous model exhibits bistability, that is, there are exactly two attractive hyperbolic fixed points---a phenomenon referred to as the Dynamical Stark Shift; and that the condition $c\in(0,4)$ is necessary and sufficient to have uniform stability for all the parameter values.

A generalized version of the model, which may be useful in certain cases, considers that the input varies over time.\cite{gammaitoni1998,bartussek1994}
For this reason we substitute the parameter $\lb$ representing the amplitude of the incident field by the parametric function $\lb+y(t)$, focusing on the cases $\lb+y\ge 0$; i.e., $\lb\geq-\inf_{r\in\R}y(r)$. We get
\begin{equation}\label{eq:parametricnonautonomousbonifacio}
x'=\lb+y(t)-x-2\,c\,\frac{x}{1+x^2}\:,
\end{equation}
which we rewrite as $x'=\lb+y(t)+g(x)$ for
\begin{equation}\label{eq:definition_g}
g(x):=-x-2\,c\,\frac{x}{1+x^2}\:.
\end{equation}
In the spirit of Refs. \onlinecite{bonifaciolugiato1,bonifaciolugiato2}, although using different methods, we examine the occurrence of uniform stability and bistability in the nonautonomous model \eqref{eq:parametricnonautonomousbonifacio} as $\lb$ varies.
In this nonautonomous scalar case, we define these terms using uniformly separated attractive hyperbolic solutions:
we say that two solutions $b_1(t)$ and $b_2(t)$ of \eqref{eq:parametricnonautonomousbonifacio} are {\em uniformly separated} if they are bounded and $\inf_{t\in\R}|b_1(t)-b_2(t)|>0$, and that a bounded solution $b(t)$ is {\em hyperbolic attractive} (resp.~{\em hyperbolic repulsive}) if there exist $k\geq1$ and $\gamma>0$ such that $\exp(\int_s^t g'(b(r))\, dr)\leq k\,e^{-\gamma(t-s)}$ for $t\ge s$ (resp.~$\exp(\int_s^t g'(b(r))\, dr)\leq k\,e^{\gamma(t-s)}$ for $t\le s$). This leads to the following definitions:
\begin{defi}
Eq.~\eqref{eq:parametricnonautonomousbonifacio}$_\lb$ is said to exhibit  {\em uniform stability} if there exists exactly one bounded solution and it is hyperbolic attractive,
and it is said to exhibit {\em (uniform) bistability} if there exist exactly two attractive hyperbolic solutions which are uniformly separated and positive.
\end{defi}

In general, a dynamical system is called bistable if it possesses two local attractors---typically, in our scalar setting, two attractive hyperbolic solutions---which may be separated by various types of dynamical structures.\cite{kuznetsov2012} In all the bistability cases in this paper, the two attractive hyperbolic solutions are separated by a repulsive hyperbolic solution, which serves as a boundary between their respective basins of attraction.

Since the equations describing the out-of-equilibrium dynamics were derived using adiabatic simplifications of several physical variables, the physical validity of the nonautonomous model \eqref{eq:parametricnonautonomousbonifacio} is possibly limited to cases where the temporal variation of $y(t)$ is sufficiently slow.

If we work with parametric values $\lb\geq-\inf_{r\in\R}y(r)$, we have $\lb+y(t)+g(0)\geq0$ for all $t\in\R$, which ensures that the half-plane $x\geq0$ is positively invariant under the dynamics of the process induced by \eqref{eq:parametricnonautonomousbonifacio}.
To use the compactification arguments underlying the results of our previous works, it is necessary for $y(t)$ to be a bounded and uniformly continuous function---a hypothesis that we will assume from this point onward without further mention: as we will explain later in this section this allows us to employ a skewproduct flow formulation.\medskip

\textbf{Properties of the map $g$.} Throughout the paper, is important to keep in mind the dependence of $g$ on $c$, not reflected in the notation for simplicity. The properties of this map are key to apply arguments of previous works, part of which also require a suitable reformulation of our nonautonomous bifurcation problem.

Some of the most important points to have in mind are the concavity, convexity and d-concavity properties of the map $x\mapsto y(t)+g(x)$ for each $t\in\R$---or equivalently those of $g(x)$---, which we analyze by computing the first derivatives of $g$:
\begin{equation}\label{eq:bonifacioderivatives}
\begin{split}
g'(x)&=\frac{-1-2\,c+2\,(c-1)\,x^2-x^4}{(1+x^2)^2}\,,\\[1ex]
g''(x)&=-4\,c\:\frac{x\,(x^2-3)}{(1+x^2)^3}\,,\\[1ex]
g'''(x)&=12\,c\:\frac{(x^2-3-2\sqrt{2}\,)(x^2-3+2\sqrt{2}\,)}{(1+x^2)^4}\,.
\end{split}
\end{equation}
An examination of \eqref{eq:bonifacioderivatives} reveals the following facts:
\begin{itemize}
\item $g$ is concave-convex on the positive half-plane $x\geq0$: it is strictly concave on $[\sqrt{3},\infty)$ and strictly convex on $[0,\sqrt{3}]$,
\item  $g$ is not globally d-concave (that is, $g'$ is not globally concave) on $[0,\infty)$, but there exists a d-concavity band around the graph of the inflection curve $x=\sqrt{3}$, independent of $c$: $g'$ is strictly concave on
\[
\left[\sqrt{3-2\sqrt{2}}\,,\:\sqrt{3+2\sqrt{2}}\,\right]\subset\left(0\,,\,\sqrt{6}\,\right)\,.\]
\end{itemize}
The existence, properties, and implications of d-concavity bands around the curve where concavity transitions to convexity are discussed---in the skewproduct framework---in Section 4.1 of Ref.~\onlinecite{dno6}.

Due to our previous studies on concave-convex nonlinearities (see Ref.~\onlinecite{dno6}), from this point forward, we will work with
\begin{equation}\label{eq:definition_barg}
\bar g(x):=\left\{
\begin{array}{ll}
g(x)&\text{if }x\geq0\,,\\[1ex]
-(1+2\,c)\,x-x^3&\text{if }x<0
\end{array}
\right.
\end{equation}
instead of $g$ (see \eqref{eq:definition_g}): $g$ is extended for $x<0$ in the simplest way so that the extension $\bar g$ is $C^2$, strictly concave in $[\sqrt3,\infty)$ and strictly convex in $(-\infty,\sqrt3]$.
The modified model is
\begin{equation}\label{eq:parametricnonautonomous_bar_g}
x'=\lb+y(t)+\bar g(x)\,,
\end{equation}
and it is important to emphasize that for $x\geq0$---the region where the solutions hold physical significance that, as previously mentioned, is positively invariant---its dynamics is identical to that of \eqref{eq:parametricnonautonomousbonifacio}.\medskip

\textbf{Skewproduct flow formulation.} Now, to equip ourselves with a richer set of dynamical tools and methods, we include \eqref{eq:parametricnonautonomous_bar_g} in a family of nonautonomous equations which allow the definition of a flow: one equation for each element $\w$ of the hull $\W_y$ of the map $y$, given by
\begin{equation}\label{eq:hull_definition}
\W_y=\mathrm{closure}\{y{\cdot}s\mid\,s\in\R\}\,,
\end{equation}
where the closure is taken in the compact-open topology and $y{\cdot}s$ denote the time-shifts of the function $y(t)$: $t\mapsto y{\cdot}s(t):=y(t+s)$ for $s\in\R$.
The conditions on boundedness and uniform continuity of $y(t)$ imply that $\W_y$ is a compact metric space endowed with a continuous flow $\sigma\colon\R\times\W_y\to\W_y$, $(s,\w)\mapsto \w{\cdot}s$ (see Theorem IV.3 of Ref.~\onlinecite{sell1971}).
So, we get the family
\begin{equation}\label{eq:parametric_nonautonomous_hull_bonifacio}
x'=\lb+\w(t)+\bar g(x)\,,\quad\w\in\W_y\,.
\end{equation}
The (possibly local) skewproduct flow is
\[
\tau_\lb(t,\w,x):=(\w{\cdot}t,u_\lb(t,\w,x))\,,
\]
where $u_\lb(t,\w,x)$ is the solution of \eqref{eq:parametric_nonautonomous_hull_bonifacio}$_\lb$ determined by $u_\lb(0,\w,x)=x$.
Note that \eqref{eq:parametric_nonautonomous_hull_bonifacio} for $\w=y\in\W_y$ is \eqref{eq:parametricnonautonomous_bar_g},
with the same dynamics on $x\ge 0$ as \eqref{eq:parametricnonautonomousbonifacio}, and that
$\lb\geq-\inf_{r\in\R}y(r)$ if and only if $\lb\ge-\inf_{r\in\R}\w(r)$ for all $\w\in\W_y$, as easily deduced from
the existence, for each $\w\in\W_y$, of a sequence $(t_n)$ such that $\w=\lim_{n\to\infty}y{\cdot}t_n$ uniformly on
compact sets.

In this framework, we say that the graph of a continuous function $\mb\colon\W_y\to\R$ is a $\tau_\lb$-\emph{copy of the base} if $\mb(\wt)=u_\lb(t,\w,\mb(\w))$ for all $\w\in\W_y$ and $t\in\R$.
The prefix $\tau_\lb$ will be omitted if there is no risk of confusion.
We say that a $\tau_\lb$-copy of the base is \emph{hyperbolic attractive} (resp. \emph{hyperbolic repulsive}) if there exist $\rho>0$, $k\geq1$ and $\gamma>0$ such that if $\w\in\W_y$ and $|\mb(\w)-x|<\rho$, then $u_\lb(t,\w,x)$ is defined for all $t\geq0$ (resp. $t\leq0$) and
$|\mb(\wt)-u_\lb(t,\w,x)|< k\,e^{-\gamma\,t}|\mb(\w)-x|$ for all $t\ge 0$ (resp. $|\mb(\wt)-u_\lb(t,\w,x)|< k\,e^{\gamma\,t}|\mb(\w)-x|$ for all $t\le 0$), which means that the graph of $\mb$ is uniformly exponentially fiber-stable at $+\infty$ (resp. $-\infty$).

\section{General results on uniform stability and bistability}\label{sec:bistability_autonomous}
All the results of this section are related to the positive critical points of the map $g$ (i.e., the unique critical points of $\bar g$, given by the zeros of $\bar g'$),
\begin{equation}\label{eq:definition_x1_and_x2}
\begin{split}
x_1(c)&:=\sqrt{c-1-\sqrt{c(c-4)}}\,,\\[1ex]
x_2(c)&:=\sqrt{c-1+\sqrt{c(c-4)}}\,,
\end{split}
\end{equation}
which are real and strictly positive if $c\ge 4$, and the associated local minimum and maximum values of $-\bar g(x)$,
\begin{equation}\label{eq:definition_lb1_and_lb2}
\begin{split}
\lb_1(c)&:=-\bar g(x_2(c))=\sqrt{\frac{c^2+10\,c-2-\sqrt{c(c-4)^3}}{2}}\:,\\
\lb_2(c)&:=-\bar g(x_1(c))=\sqrt{\frac{c^2+10\,c-2+\sqrt{c(c-4)^3}}{2}}\:,
\end{split}
\end{equation}
also real and strictly positive if $c\ge 4$.
For further purposes, we note that the functions $-x_1$, $x_2$, $\lambda_1$, and $\lambda_2$ are strictly increasing on $[4, \infty)$. This can be verified analytically with some effort and can also be readily observed by plotting their graphs using any appropriate software.

The first result demonstrates the absence of bistability of \eqref{eq:parametricnonautonomousbonifacio}$_\lb$ for every $\lb$ if $0<c\le 4$, and it is based on a property of Lyapunov exponents which is also shared by the autonomous problem \eqref{eq:autonomousbonifacio}.

\begin{prop}\label{prop:cleq4_no_bistability} If $0<c<4$, then \eqref{eq:parametricnonautonomousbonifacio}$_\lb$ exhibits uniform stability for any value of $\lb$.
In addition, the hyperbolic solution is positive if $\lb\ge-\inf_{r\in\R} y(r)$.
Finally, if $c=4$, \eqref{eq:parametricnonautonomousbonifacio}$_\lb$ does not have two uniformly separated hyperbolic solutions for any value of $\lb$.
\end{prop}
\begin{proof}
The main argument of the proof is based on the fact that $g'(x)<0$ for all $x\in\R$ if $0<c<4$:
the four zeroes of $g$ are $\pm x_1(c)$ and $\pm x_2(c)$, with $x_1(c)$ and $x_2(c)$ given
in \eqref{eq:definition_x1_and_x2}. So, they are non-real if $0<c<4$, and $g'(0)=-1-2\,c<0$.

Let us fix any $\lb$. Since
\begin{equation}\label{eq:auxiliar}
\lim_{x\to\pm\infty}\big(\lb+\w(0)+\bar g(x)\big)=\mp\infty
\end{equation}
uniformly on $\W_y$, a global attractor $\mA_\lb\subseteq\W_y\times\R$ exists for the skewproduct flow defined on $\W_y\times\R$ from \eqref{eq:parametricnonautonomousbonifacio} (see Theorem 5.1 of Ref.~\onlinecite{dno1}).
Since each Lyapunov exponent of $\mA_\lb$ is given by the integral with respect to an ergodic measure on $\W_y$ of $g'$ evaluated over a certain function with graph in $\mA_\lb$ (see Theorem 1.8.4 of Ref.~\onlinecite{arnold1998} and see Theorem 1.36 of Ref.~\onlinecite{duenastesis}), all these exponents are strictly negative.
Consequently (see Theorem 2.13(vii) of Ref.~\onlinecite{duenastesis}, based on Theorem 3.4 of Ref.~\onlinecite{campos}), the attractor is an attractive hyperbolic copy of the base.
This means the existence of a unique (attractive hyperbolic) bounded solution for \eqref{eq:parametricnonautonomousbonifacio}$_\lb$, as asserted.

The positiveness of the hyperbolic solution when $\lb\ge-\inf_{r\in\R} y(r)$ is based on the positive invariance of the compact set
$\W_y\times[0,r]$ for a constant $r$ large enough to guarantee $\lb+\w(0)+g(r)<0$ for all $\w\in\W_y$: the attractor $\mA_\lb$ is contained in this set (see, e.g.,
Theorem 2.2 of Ref.~\onlinecite{cheban2002}).

Finally, if two separated hyperbolic solutions were to exist for $c=4$, then the persistence of hyperbolic solutions of \eqref{eq:parametricnonautonomousbonifacio} under small parametric perturbations (which is a classical result, see, e.g., Theorem 3.8 of Ref.~\onlinecite{potzsche2011}) would guarantee it also for $c<4$ close enough, which is not possible.
\end{proof}

The second result establishes a first condition on the function $y(t)$ which ensures the existence of bistability for certain parameter values, provided that $c>4$.
This result is derived using comparison arguments with the equilibria of the autonomous problem \eqref{eq:autonomousbonifacio} combined with Lyapunov exponent arguments. Recall that $\lb_1(c)>0$ if $c\ge 4$, and hence, if $\lb\ge \lb_1(c)-\inf_{r\in\R}y(r)$, then the invariance condition $\lb>-\inf_{r\in\R} y(r)$ holds. We define
\begin{equation}\label{def:h1}
 h_1(c):=\lb_2(c)-\lb_1(c)=-\bar g(x_1(c))+\bar g(x_2(c))\,,
\end{equation}
and observe that $h_1(c)>0$ for $c>4$: see definitions \eqref{eq:definition_lb1_and_lb2}.

\begin{teor}\label{teor:bistability_comparison_autonomous}
Let $c>4$ be fixed. If
\begin{equation}\label{eq:condition_bound_y1}
\lb>\lb_1(c)-\inf_{r\in\R}y(r)\,,
\end{equation}
then there exists exactly one positive bounded solution $u_\lb$ of \eqref{eq:parametricnonautonomousbonifacio}$_{\lb}$ strictly above $x_2(c)$, it is hyperbolic attractive, and it satisfies $\lim_{t\to\infty}(u_\lb(t)-x_\lb(t))=0$ for any solution $x_\lb$ taking any value in $[x_2(c),\,\infty)$; and if
\begin{equation}\label{eq:condition_bound_y2}
\lb\in\left[-\inf_{r\in\R}y(r)\,,\;\lb_2(c)-\sup_{r\in\R}y(r)\right)\,,
\end{equation}
then there exists exactly one positive bounded solution $l_\lb$ of \eqref{eq:parametricnonautonomousbonifacio}$_{\lb}$ strictly below $x_1(c)$, it is hyperbolic attractive, it is strictly positive if $\lb>-\inf_{r\in\R}y(r)$, and it satisfies $\lim_{t\to\infty}(l_\lb(t)-x_\lb(t))=0$ for any solution $x_\lb$ taking any value in $[0,\,x_1(c)]$. In addition, if
\begin{equation}\label{cond:h1}
\sup_{r\in\R}y(r)-\inf_{r\in\R}y(r)<h_1(c)
\end{equation}
and
\begin{equation}\label{eq:lambda_interval_autonomous_prop}
\lb\in\mI_1(c):=\left(\lb_1(c)-\inf_{r\in\R}y(r)\,,\;\lb_2(c)-\sup_{r\in\R}y(r)\right)\,,
\end{equation}
then the equation \eqref{eq:parametricnonautonomousbonifacio}$_{\lb}$ exhibits bistability.
More precisely, the basins of attraction on the positive half-plane of the two previously found attractive hyperbolic solutions are separated by the graph of a repulsive hyperbolic solution $m_\lb$ which is above $x_1(c)$ and below $x_2(c)$.

Finally, the last conclusions hold if $\sup_{r\in\R}y(r)-\inf_{r\in\R}y(r)\le h_1(c)$,
$\lb\in\bar\mI_1(c):=\big[\,\lb_1(c)-\inf_{r\in\R}y(r),$ $\lb_2(c)-\sup_{r\in\R}y(r)\,\big]$, and the constant maps $\inf_{r\in\R}y(r)$ and $\sup_{r\in\R}y(r)$ are not in $\W_y$.
\end{teor}
\begin{proof}
The proof relies on the properties of the bifurcation diagram for the autonomous model $x'=\lb+\bar g(x)$.
The fixed points of \eqref{eq:autonomousbonifacio} for each $\lb\geq0$ are the solutions of $\lb=-g(x)$,
which, for $x\geq0$, describe an S-shaped curve as a function of $\lb\geq0$: the bifurcation diagram is depicted in Fig.~\ref{fig:bif_diagram_aut_bonifacio}.
 \begin{figure}
     \centering
         \includegraphics[width=\columnwidth]{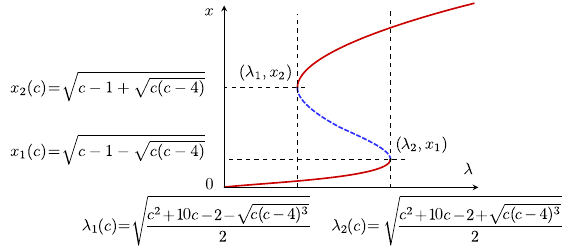}
         \caption{Bifurcation diagram for the autonomous model \eqref{eq:autonomousbonifacio} for $x\geq0$
         and any fixed $c>4$: in this drawing, $c=5$.
         The hyperbolic attractive fixed points are represented in red, while the hyperbolic repulsive fixed points are shown in blue.
         The non-hyperbolic fixed points, which correspond to bifurcation points, are depicted in black.}
        \label{fig:bif_diagram_aut_bonifacio}
\end{figure}
The values $\lb_1(c)$ and $\lb_2(c)$ of the parameter are the unique (saddle-node) bifurcation points in the autonomous case.

Now we will use comparison arguments.
First, since the hull construction in Section \ref{sec:bonifacio-lugiato} guarantees that $\inf_{r\in\R}y(r)\leq\w(0)\leq\sup_{r\in\R}y(r)$
for all $\w\in\W_y$, and since $\lb_1(c)=-\bar g(x_2(c))$, condition \eqref{eq:condition_bound_y1} ensures that
\begin{equation}\label{eq:x2_strict_lower_solution}
0<\lb+\w(0)+\bar g(x_2(c))
\end{equation}
for all $\w\in\W_y$.
That is, the constant $x_2(c)$ is a strict global lower solution for \eqref{eq:parametric_nonautonomous_hull_bonifacio}$_{\lb}$
(see Section 2.2 of Ref.~\onlinecite{dno1}).
Since $\lim_{x\to\infty}\bar g(x)=-\infty$, there exists a large enough constant $r_0>x_2(c)$ such that every $r\geq r_0$ is a strict global upper solution of \eqref{eq:parametric_nonautonomous_hull_bonifacio}$_{\lb}$.
So, for any $r\geq r_0$, $\W_y\times[x_2(c),r]$ is a compact forward invariant set for the corresponding skewproduct flow $\tau_{\lb}$ on $\W_y\times\R$.
This ensures the existence of a local attractor $\mA^u_{\lb}$ contained in $\W_y\times[x_2(c),r]$ (see Theorem 2.2 of Ref.~\onlinecite{cheban2002}), which attracts all the $\tau_{\lb}$-orbits taking any value at $\W_y\times[x_2(c),r]$, and which in addition is composed by all the globally defined bounded $\tau_{\lb}$-orbits which lie within $\W_y\times[x_2(c),r]$ (see Lemma 1.6 of Ref.~\onlinecite{carvalho2013}).
It is easy to deduce from the strict character of the global lower and upper solutions $x_2(c)$ and $r$ that this local attractor is contained on $\W_y\times(x_2(c),r)$.
Hence, since $\bar g'(x)<0$ on this region, the upper Lyapunov exponent of $\mA^u_{\lb}$ is strictly negative (see again Theorem 1.36 of Ref.~\onlinecite{duenastesis}).
By reasoning as in Theorem 3.4 of Ref.~\onlinecite{campos} and Theorem 2.13(vii) of Ref.~\onlinecite{duenastesis}, we conclude that $\mA^u_{\lb}$ is an attractive hyperbolic copy of the base.
Therefore, $\mA^u_{\lb}$ is independent of the choice of $r\geq r_0$.
Hence, the unique bounded solution $u_\lb$ for \eqref{eq:parametricnonautonomousbonifacio}$_\lb$ above $x_2(c)$ is an attractive hyperbolic solution, it is placed strictly above $x_2(c)$, and it attracts as time increases any solution of \eqref{eq:parametricnonautonomousbonifacio}$_\lb$ with initial data greater or equal than $x_2(c)$ at any initial time, as asserted.

An analogous argument shows the claim when we assume \eqref{eq:condition_bound_y2}: transforming $\lb<\lb_2(c)-\sup_{r\in\R}y(r)$ into
\begin{equation}\label{eq:x1_strict_upper_solution}
0>\lb+\w(0)+\bar g(x_1(c))
\end{equation}
for all $\w\in\W_y$, we deduce that in this case $\W_y\times[0,\,x_1(c)]$ is a compact forward invariant set for $\tau_{\lb}$, and the existence of a local attractor $\mA^l_{\lb}$ contained in $\W_y\times[0,\,x_1(c))$ (or in $\W_y\times(0,\,x_1(c))$ if $\lb>-\inf_{r\in\R}y(r)$) which is an attractive hyperbolic copy of the base.

It is obvious that the interval $\mI_1(c)$ given by \eqref{eq:lambda_interval_autonomous_prop} is non-degenerate if \eqref{cond:h1} holds. Clearly, every $\lambda\in\mI_1(c)$ fulfills both \eqref{eq:condition_bound_y1} and \eqref{eq:condition_bound_y2}, implying the existence of two strictly positive attractive hyperbolic solutions of \eqref{eq:parametric_nonautonomous_hull_bonifacio}$_\lb$: $u_\lb$ above $x_2(c)$ and $l_\lb$ below $x_1(c)$.
In addition, conditions \eqref{eq:x2_strict_lower_solution} and \eqref{eq:x1_strict_upper_solution} ensure that $\W_y\times[x_1(c),\,x_2(c)]$ is a compact backward invariant set for $\tau_{\lb}$.
Therefore, $\W_y\times[x_1(c),\,x_2(c)]$ is forward invariant for the time-reversed flow of $\tau_{\lb}$
(see the definition in Section 2.2 of Ref.~\onlinecite{dno1}).
Thus, the reasoning of the previous part of the proof shows that there exists a local attractor $\mA^m_{\lb}$ for the time-reversed flow of $\tau_{\lb}$ contained in $\W_y\times(x_1(c),\,x_2(c))$.
Since the time-reversed family of equations of \eqref{eq:parametric_nonautonomous_hull_bonifacio}$_{\lb}$ is
\[
\frac{dx}{ds}=-\lb-\w(-s)-\bar g(x)\,,\quad\w\in\W_y\,,
\]
where $s=-t$, and $-\bar g'(x)<0$ for $x\in(x_1(c),\,x_2(c))$, the upper (time-reversed) Lyapunov exponent of $\mA^m_{\lb}$ is strictly negative.
Therefore, $\mA^m_{\lb}$ is a (time-reversed) attractive hyperbolic copy of the base, so it is a repulsive hyperbolic copy of the base for $\tau_{\lb}$, and it is composed by all the bounded orbits contained in $\W_y\times[x_1(c),\,x_2(c)]$.
Hence, there is only one bounded solution $m_\lb$ of \eqref{eq:parametricnonautonomousbonifacio}$_\lb$ taking all its values in $[x_1(c),\,x_2(c)]$, it is hyperbolic repulsive, and it attracts every solution taking any value in $[x_1(c),\,x_2(c)]$ as time decreases.
Since any solution taking the value $x_2(c)$ (resp.~$x_1(c)$) approaches $u_\lb$ (resp.~$l_\lb$) as time increases and $m_\lb$ as time decreases, it cannot be hyperbolic (see Lemma 3.8 of Ref.~\onlinecite{dno3}).
So, there are exactly two positive attractive hyperbolic solutions; i.e., there is bistability.
On the other hand, our proof shows that any bounded solution strictly above (resp. below) $m_\lb$ and below $u_\lb$ (resp.~above $l_\lb$) takes the value $x_2(c)$ (resp. $x_1(c)$) at a certain time, and hence it approaches $u_\lb$ (resp.~$l_\lb$) as time increases.

Let us prove the last assertion of the theorem. If $\lb$ satisfies \eqref{eq:lambda_interval_autonomous_prop}, there is nothing to prove. Let us work with $\lb_0:=\lb_1(c)-\inf_{r\in\R}y(r)$. Reasoning as at the beginning of this proof, we establish the existence of a local attractor $\mA^u_{\lb_0}$ in $\W_y\times[x_2(c),r]$ for an $r$ large enough. To prove that the upper Lyapunov exponent of $\mA^u_{\lb_0}$ is strictly negative, we must check that $\int_{\W_y}\bar g'(\mb(\w))\,dm<0$ for any ergodic measure $m$ on $\W_y$ and any $m$-measurable map $\mb\colon\W_y\to\R$ with graph contained in $\mA^u_{\lb_0}$ such that $t\mapsto\mb(\wt)$ is $C^1$ and solves \eqref{eq:parametric_nonautonomous_hull_bonifacio}$_{\lb_0}$ for all $\w\in\W_y$ (see, e.g., Theorem 1.36 of Ref.~\onlinecite{duenastesis}). Recall that $\bar g'(x)\le 0$ for all $(\w,x)\in \mA^u_{\lb_0}$. Let us fix $\w_0$ in the support of $m$. It is not hard to check that, if $\mb(\w_0{\cdot}t)=x_2(c)$ for some $t\in\R$, then $\w_0(t)=\inf_{r\in\R}y(r)$. Since the constant map $\inf_{r\in\R}y(r)$ does not belong to $\W_y$, there exists $t_0$ such that $\mb(\w_0{\cdot}t_0)>x_2(c)$, which combined with the continuity of $\bar g'$ ensures that $\bar g'(\w)\le-\delta<0$ for all the points $\w$ in an open ball $\mB$ of $\W_y$. The set $\mB_{-t_0}:=\{\sigma(-t_0,\w)\mid\w\in\mB\}\ni\w_0$ is open, and hence $m(\mB)=m(\mB_{-t_0})>0$. So, $\int_{\W_y}\bar g'(\mb(\w))\,dm\le \int_{\mB}\bar g'(\mb(\w))\,dm\le-\delta\,m(\mB)<0$, as asserted. The rest of the proof for $\lb_0$ uses similar arguments and those used before in this proof, and the arguments for $\lb_2(c)-\sup_{r\in\R}y(r)$ are analogous.
\end{proof}
\begin{remark}
\label{rm:nueva}
1.~Note that, although the input $I(t):=\lb+y(t)$ does not uniquely determine $\lb$ and $y(t)$, the difference $\sup_{r\in\R}(\lb+y(r))-\inf_{r\in\R}(\lb+y(r))$ is uniquely determined.

2.~Condition \eqref{eq:lambda_interval_autonomous_prop} is equivalent to say that the input map $\lb+y(t)$ takes values in a compact subset of $(\lb_1(c),\,\lb_2(c))$.
In fact, Theorem~\ref{teor:bistability_comparison_autonomous} can be directly applied to equation \eqref{eq:intro} to guarantee that, if the input $I(t)$ takes values in a compact subset of the interval $(\lb_1(c),\,\lb_2(c))$ determining the bistability interval for \eqref{eq:autonomousbonifacio}, then \eqref{eq:intro} exhibits bistability.

3.~Arguments analogous to those of Theorem \ref{teor:bistability_comparison_autonomous} allow us to prove that, if $\lb+y(t)$ takes values in a compact subset of $[0,\,\lb_1(c))$, then \eqref{eq:parametricnonautonomousbonifacio}$_\lb$ exhibits uniform stability with a unique positive attractive hyperbolic solution that is strictly below $x_1(c)$; and that, if $\lb+y(t)$ takes values in a compact subset of $(\lb_2(c),\infty)$, then \eqref{eq:parametricnonautonomousbonifacio}$_\lb$ exhibits uniform stability with a unique (strictly positive) attractive hyperbolic solution that is strictly above $x_2(c)$. As in the previous remark, these results can be extended to \eqref{eq:intro} in terms of the values taken by $I(t)$.
\end{remark}

\section{Further results on bistability: auxiliary concave-linear and linear-convex equations}\label{sec:bistability_concaveconvex}
In this section, we draw on some of the tools developed in Ref.~\onlinecite{dno6} to broaden the results presented in Section~\ref{sec:bistability_autonomous}, albeit under significantly more restrictive assumptions regarding the variation of $y(t)$.
Specifically, we get information for our concave-convex equation from the bifurcation diagrams of two auxiliary concave-linear and linear-convex equations. To formulate the auxiliary equations as described in Ref.~\onlinecite{dno6}, first, we shift the inflection curve $x=\sqrt{3}$ in \eqref{eq:parametric_nonautonomous_hull_bonifacio} (common for all $\w\in\W_y$) to $0$ by means of the change of variables $z=x-\sqrt{3}\,$:
\begin{equation}\label{eq:change_of_variables_x_z}
z'=\lb+\w(t)+\bar g(z+\sqrt{3}\,)\,.
\end{equation}
So, the region $x\geq 0$ we are interested in is taken to $z\geq-\sqrt{3}$.
And second, we define $\mg$ so as to get
\begin{equation}\label{eq:construction_of_mg}
\begin{split}
\bar g(z+\sqrt{3}\,)
&=\bar g(\sqrt{3})+\tilde g'(\sqrt{3})\,z+\mg(z)\\[1ex]
&=-\frac{\sqrt{3}\,(c+2)}{2}+\left(\frac{c}{4}-1\right)\,z+\mg(z)
\end{split}
\end{equation}
and rewrite \eqref{eq:change_of_variables_x_z} as
\begin{equation}\label{eq:concave-convex_bonifacio}
z'=\lb+\mc(\wt)+\md\,z+\mg(z)\,,
\end{equation}
so that
\begin{equation}\label{eq:definition_mg}
\begin{split}
\mc(\w)&:=\w(0)-\frac{\sqrt{3}\,(c+2)}{2}\:,\\
\md&:=c/4-1\,,\\
\mg(z)&:=\bar g(z+\sqrt{3}\,)+\frac{\sqrt{3}\,(c+2)}{2}-\md\,z\,.
\end{split}
\end{equation}
The family \eqref{eq:concave-convex_bonifacio}
also induces a skewproduct flow. In addition,
\begin{lema}
\label{lema:cc1-cc6}
If $c>4$\,, then
\begin{itemize}
\item $\md>0$, $\mc\in C(\W,\R)$, $\mg\in C^2(\R,\R)$\,,
\item $\lim_{z\to\pm\infty}(\md\,z+\mg(z))=\mp\infty$\,,
\item $\mg(0)=\mg'(0)=\mg''(0)=0$\,,
\item $\mg(z)<0$ for all $z>0$ and $\mg(z)>0$ for all $z<0$\,,
\item $z\mapsto \mg(z)$ is strictly concave on $[0,\infty)$ and strictly convex on $(-\infty,0]$\,.
\end{itemize}
\end{lema}
\begin{proof}$\md>0$ follows directly from $c>4$, the continuity of $\mc$ is clear, and $\mg\in C^2(\R,\R)$ follows from \eqref{eq:definition_mg} and $\bar g\in C^2(\R,\R)$.
In addition,
\[
\lim_{z\to\pm\infty}(\md\,z+\mg(z))=\lim_{z\to\pm\infty}\tilde g(z+\sqrt{3})+\frac{\sqrt{3}\,(c+2)}{2}=\mp\infty\,,\]
and \eqref{eq:construction_of_mg} ensures $\mg(0)=\mg'(0)=0$ and $\mg''(0)=\bar g''(\sqrt{3})=g''(\sqrt{3})=0$ (see \eqref{eq:bonifacioderivatives}).

Since $\mg''(z)=\bar g''(z+\sqrt3)$, the strict concavity and convexity properties of $\mg$ at $z\ge 0$ and $z\le 0$ follow from those of $\bar g$ on $x\ge\sqrt3$ and $x\le\sqrt3$.
Thus, $\mg'$ is strictly decreasing on $[0,\infty)$ and strictly increasing on $(-\infty,0]$, which
combined with $\mg(0)=\mg'(0)=0$ yields $\mg(z)<0$ for all $z>0$ and $\mg(z)>0$ for all $z<0$.
\end{proof}
\begin{remark}Lemma \ref{lema:cc1-cc6} shows that, if $c>4$, then the map $\mf(\w,z):=\mc(\w) + \md\,z+\mg(z)$ satisfies the hypotheses \textbf{cc1-cc6} described at the beginning of Section 4 of Ref.~\onlinecite{dno6}.
\end{remark}

We now introduce the concave-linear and linear-convex equations, which are obtained by replacing the function $\mg$ in equation \eqref{eq:concave-convex_bonifacio} with respective functions $\mg_-$ and $\mg_+$.
These functions respectively vanish in the half-planes $z<0$ and $z>0$ and coincide with $\mg$ where they do not vanish. That is,
\[
\mg_-(z):=\left\{\begin{array}{ll}
\mg(z)
&\text{if }z\geq0\,,\\[1.5ex]
0&\text{if }z<0\,,
\end{array}\right.
\quad\;
\mg_+(z):=\left\{\begin{array}{ll}
0&\text{if }z>0\,,\\[1.5ex]
\mg(z)
&\text{if }z\leq0\,.
\end{array}\right.
\]
The concave-linear family is
\begin{equation}\label{eq:concave-linear_bonifacio_bifurcation}
z'=\lb+\mc(\wt)+\md\,z+\mg_-(z)\,,\quad\w\in\W_y\,,
\end{equation}
and the linear-convex family is
\begin{equation}\label{eq:linear-convex_bonifacio_bifurcation}
z'=\lb+\mc(\wt)+\md\,z+\mg_+(z)\,,\quad\w\in\W_y\,.
\end{equation}
Theorems 3.4 and 3.5 of Ref.~\onlinecite{dno6} establish the existence of unique bifurcation point $\lb_-(c)\in\R$, such that \eqref{eq:concave-linear_bonifacio_bifurcation}$_\lb$ admits two uniformly separated hyperbolic solutions---attractive the upper one and repulsive the lower one---for all $\w\in\W_y$ if $\lb>\lb_-(c)$, whereas this is not the case for $\w=y\in\W_y$ if $\lb\leq\lb_-(c)$.
On the other hand, Remark 3.6 of Ref.~\onlinecite{dno6} shows the symmetrical properties for the linear-convex bifurcation problem \eqref{eq:linear-convex_bonifacio_bifurcation}$_\lb$:
the existence of a uniquely defined bifurcation point $\lb_+(c)$ such that \eqref{eq:linear-convex_bonifacio_bifurcation}$_\lb$ admits two uniformly separated hyperbolic solutions---attractive the lower one and repulsive the upper one---for all $\w\in\W_y$ if $\lb<\lb_+(c)$, whereas this is not the case for $\w=y\in\W_y$ if $\lb\geq\lb_+(c)$.

\begin{lema}\label{lema:intervalo}
Let $c>4$ be fixed. Then, $\lb_-(c)\ge\lb_1(c)-\sup_{r\in\R}y(r)$ and $\lb_+(c)\le\lb_2(c)-\inf_{r\in\R}y(r)$.
\end{lema}
\begin{proof}
The change of variables $x:=z+\sqrt3$ takes \eqref{eq:concave-linear_bonifacio_bifurcation} to $x'=\lb+\w(t)+\bar g_-(x)$, where $\bar g_-(x)$ coincides with $\bar g$ on $[\sqrt3,\infty)$ and takes the value $\bar g(\sqrt3)+\bar g'(\sqrt3)(x-\sqrt3)$ on $(-\infty,\sqrt3)$. The global maximum of $\bar g_-$, reached at $x_2(c)$, is $-\lb_1(c)$. So, if $\lb<\lb_1(c)-\sup_{r\in\R}y(r)$, then there exists $\delta>0$
such that every solution of $x'=\lb+y(t)+\bar g_-(x)$ satisfies $x'(t)<-\delta$ for all $t\in\R$, and this precludes the existence of bounded solutions. It follows that $\lb_-(c)\ge\lb_1(c)-\sup_{r\in\R}y(r)$. The second  inequality is proved similarly.
\end{proof}

Two additional relevant parameters associated with these bifurcation problems are $\mu_-(c)$ and $\mu_+(c)$.
The parameter $\mu_-(c)$ represents the infimum of the values of $\lb$ for which the unique bounded solution of the linear problem
\[
 z'=\lb+\mc(\wt)+\md\,z
\]
---which is hyperbolic repulsive since $\md>0$, and decreases with $\lb$---is non-positive for all $\w\in\W_y$.
Note that this solution is $t\mapsto b^\lb_\w(t):=\mb^\lb(\w{\cdot}t)$, where
\[
\mb^\lb(\w):=\frac{1}{\md}\left(\frac{\sqrt{3}\,(c+2)}{2}-\lb\right)-\int_0^\infty e^{-\md\,s}\w(s)\,ds
\]
(see Chapter 3 of Ref.~\onlinecite{coppel1978}).
In turn, $\mu_+(c)\leq\mu_-(c)$ represents the supremum of the values of $\lb$ for which this solution is non-negative for all $\w\in\W_y$.
Therefore,
\begin{equation}\label{eq:mu-_mu+_definition}
\begin{split}
\mu_-(c)&=\frac{\sqrt{3}\,(c+2)}{2}-\md\inf_{\w\in\W_y}\int_0^\infty e^{-\md\,s}\w(s)\,ds\,,\\
\mu_+(c)&=\frac{\sqrt{3}\,(c+2)}{2}-\md\sup_{\w\in\W_y}\int_0^\infty e^{-\md\,s}\w(s)\,ds\,.
\end{split}
\end{equation}

This completes the preliminaries needed to get new results on bistability for \eqref{eq:parametricnonautonomousbonifacio} based on the properties of \eqref{eq:concave-linear_bonifacio_bifurcation} and \eqref{eq:linear-convex_bonifacio_bifurcation}.

The following proposition shows that the interval $\mI_1(c)$ of \eqref{eq:lambda_interval_autonomous_prop} is contained in $(\lb_-(c),\,\lb_+(c))$ for $c>4$.
\begin{prop}\label{prop:encaje}
Let $c>4$ be fixed. Then, for all $\lb>\lb_-(c)$ the upper bounded solution of \eqref{eq:parametricnonautonomous_bar_g} is above $b_y^\lb$ and hyperbolic attractive; and for all $\lb<\lb_+(c)$ the lower bounded solution
of \eqref{eq:parametricnonautonomous_bar_g} is below $b_y^\lb$ and hyperbolic attractive. In addition, if \eqref{cond:h1} holds, then
\[
\begin{split}
\lb_1(c)-\sup_{r\in\R}y(r)\le&\,\lb_-(c)\le\lb_1(c)-\inf_{r\in\R}y(r)\,,\\
\lb_2(c)-\sup_{r\in\R}y(r)\le&\,\lb_+(c)\le\lb_2(c)-\inf_{r\in\R}y(r)\,.
\end{split}
\]
\end{prop}
\begin{proof}
Proposition 4.3 of Ref.~\onlinecite{dno6} ensures the first assertions. Let us assume \eqref{cond:h1}. The inequality $\lb_1(c)-\sup_{r\in\R}y(r)\le\lb_-(c)$ is proved in Lemma \ref{lema:intervalo}. To check that $\lb_-(c)\le\lb_1(c)-\inf_{r\in\R}y(r)$, we take $\lb>\lb_1(c)-\inf_{r\in\R}y(r)$. By reviewing the proof of Theorem~\ref{teor:bistability_comparison_autonomous}, we observe that there exists an attractive hyperbolic copy of the base for $\tau_\lb$, which in addition is above $x_2(c)>\sqrt3$. The change of variables $(\w,x)\mapsto (\w,z)$ with $z:=x-\sqrt3$ takes it to an attractive hyperbolic copy of the base for the flow induced by \eqref{eq:concave-convex_bonifacio}$_\lb$. Since this copy of the base is above $0$, it is also an attractive hyperbolic copy of the base for the flow induced by the family \eqref{eq:concave-linear_bonifacio_bifurcation}$_\lb$. Theorem 3.4 of Ref.~\onlinecite{dno6} shows that $\lb>\lb_-(c)$, from where the inequality follows. Lemma \ref{lema:intervalo} proves the first inequality of the second chain. The second one can be proved with an argument similar to that just used, having in mind that $x_1(c)<\sqrt3$.
\end{proof}

The next goal is to find conditions ensuring that, in fact, $(\lb_-(c),\,\lb_+(c))\cap[-\inf_{r\in\R}y(r),\,\lb_+(c))$ is the bistability interval for our problem. To this end, we will establish conditions (in general) more restrictive than those of Proposition~\ref{prop:encaje} ensuring that
\begin{equation}\label{eq:desmu}
\begin{split}
 &\lb_1(c)-\inf_{r\in\R}y(r)<\mu_+(c)\,,\\
 &\mu_-(c)<\lb_2(c)-\sup_{r\in\R}y(r)\,.
\end{split}
\end{equation}
The proof of Theorem \ref{teor:concave_lineal_bounds} will clarify the interest of these inequalities. According to \eqref{eq:definition_lb1_and_lb2} and \eqref{eq:mu-_mu+_definition}, and since $\bar g(\sqrt3)=-\sqrt{3}\,(c+2)/2$,
the first inequality in \eqref{eq:desmu} is equivalent to
\[
 \sup_{r\in\R}\,\md\int_0^\infty e^{-\md\,s}y(r+s)\,ds-\inf_{r\in\R}y(r)<h_2(c)
\]
for
\begin{equation}\label{def:h2}
 h_2(c):=\bar g(x_2(c))-\bar g(\sqrt3)=\sqrt{3}\,(c+2)/2-\lb_1(c)\,,
\end{equation}
and the second one to
\[
 \sup_{r\in\R}y(r)-\inf_{r\in\R}\,\md\int_0^\infty e^{-\md\,s}y(r+s)\,ds<h_3(c)
\]
for
\begin{equation}\label{def:h3}
 h_3(c):=\bar g(\sqrt3)-\bar g(x_1(c))=\lb_2(c)-\sqrt{3}\,(c+2)/2\,.
\end{equation}
The explicit expressions of $h_2(c)$ and $h_3(c)$,
\[
\begin{split}
 h_2(c)&=\frac{\sqrt{3}\,(2+c)}{2}
 -\sqrt{\frac{c^2+10\,c-2-\sqrt{c(c-4)^3}}{2}}\:,\\
 h_3(c)&=\sqrt{\frac{c^2+10\,c-2+\sqrt{c(c-4)^3}}{2}}-\frac{\sqrt{3}\,(2+c)}{2}\:,
\end{split}
\]
can be useful when trying to apply the results of this section to a particular example.
\begin{lema}\label{lema:order}
Let $h_1,\,h_2,\,h_3\colon[4,\infty)\to\R$ be the maps given by \eqref{def:h1}, \eqref{def:h2} and \eqref{def:h3}. Then, $h_1(4)=h_2(4)=h_3(4)=0$ and $h_1(c)>h_2(c)>h_3(c)>0$ for all $c>4$.
\end{lema}
\begin{proof}
The equalities for $c=4$ follow from $x_1(4)=x_2(4)=\sqrt3$. Assume now $c>4$. Since $\bar g'(x)=g'(x)>0$ on $(x_1(c),\,x_2(c))$ and $x_1(c)<\sqrt3<x_2(c)$, we get $\bar g(x_1(c))<\bar g(\sqrt3)<\bar g(x_2(c))$, i.e.,
$-\lb_2(c)<-\sqrt3(c+2)/2<-\lb_1(c)$, which yields $h_1(c)>h_2(c)$ and $h_3(c)>0$. Finally, $h_2(c)>h_3(c)$ is equivalent to $\lb_1(c)+\lb_2(c)>\sqrt3\,(c+2)$. Some standard manipulation (involving squaring twice) shows that this is equivalent to  $(c-4)^2(c^2+2c+3)>0$, true for $c>4$.
\end{proof}
\begin{figure}
     \centering
         \includegraphics[width=\columnwidth]{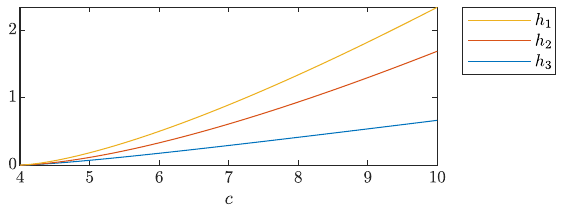}
         \caption{Representation of the functions $h_1(c)$, $h_2(c)$ and $h_3(c)$ given by \eqref{def:h1}, \eqref{def:h2} and \eqref{def:h3} on $[4,10]$.}
        \label{fig:j_functions}
\end{figure}

Figure \ref{fig:j_functions} depicts the maps $h_1,\,h_2$ and $h_3$. We can state the main result of this section, whose hypotheses, as explained below, also ensure the occurrence of one or two nonautonomous saddle-node bifurcation points of hyperbolic solutions.
\begin{teor}\label{teor:concave_lineal_bounds}
Let $c>4$ be fixed.
If
\begin{equation}\label{eq:cond_mu+}
 \sup_{r\in\R}\,\md\int_0^\infty e^{-\md\,s}y(r+s)\,ds-\inf_{r\in\R}y(r)<h_2(c)\,,
\end{equation}
with $h_2$ defined by \eqref{def:h2}, then \eqref{eq:condition_bound_y1} holds for $\lb=\mu_+(c)$, and hence there exists a strictly positive attractive hyperbolic solution $u_\lb$ of \eqref{eq:parametricnonautonomousbonifacio}$_\lb$ for all $\lb\ge\mu_+(c)$. If
\begin{equation}\label{eq:cond_mu-}
 \sup_{r\in\R}y(r)-\inf_{r\in\R}\,\md\int_0^\infty e^{-\md\,s}y(r+s)\,ds<h_3(c)\,,
\end{equation}
with $h_3$ defined by \eqref{def:h3}, then $\mu_-(c)<\lb_2(c)-\sup_{r\in\R}y(r)$, and thus \eqref{eq:condition_bound_y2} holds for $\lb=\mu_-(c)$ if $\mu_-(c)\ge-\inf_{r\in\R}y(r)$. If so,
there exists a positive attractive hyperbolic solution $l_\lb$ of \eqref{eq:parametricnonautonomousbonifacio}$_\lb$ for all $\lb\in[-\inf_{r\in\R}y(r),\,\mu_-(c)]$. And, if both \eqref{eq:cond_mu+} and \eqref{eq:cond_mu-} hold and $\mu_-(c)\ge-\inf_{r\in\R}y(r)$,
then \eqref{eq:parametricnonautonomousbonifacio}$_\lb$ exhibits bistability for $\lb\in\mI_2(c):=(\lb_-(c),\,\lb_+(c))\cap[-\inf_{r\in\R}y(r),\,\lb_+(c))$ and not for any other $\lb\geq-\inf_{r\in\R}y(r)$.
If, in addition, $\sup_{r\in\R}y(r)-\inf_{r\in\R}y(r)\le\lb_1(c)$, then $\mI_2(c)=(\lb_-(c),\,\lb_+(c))$.

More precisely: for $\lb\in\mI_2(c)$ there exist three positive hyperbolic solutions of \eqref{eq:parametricnonautonomousbonifacio}$_\lb$, $l_\lb<m_\lb<u_\lb$, with $l_\lb$ and $u_\lb$ attractive and with basins of attraction on the positive half-plane separated by the graph of $m_\lb$, which is repulsive;
if $\lb\ge\lb_+(c)$, then $u_\lb$ is the unique positive hyperbolic solution, it is attractive, and $(\lb_-(c),\,\infty)\cap[-\inf_{r\in\R}y(r),\,\infty)\to C(\R,\R),\,\lb\mapsto u_\lb$ is continuous in the uniform topology; and if $-\inf_{r\in\R}y(r)\le\lb_-(c)$ and $\lb\in[-\inf_{r\in\R}y(r),\lb_-(c)]$, then $l_\lb$ is the unique positive hyperbolic solution, it is attractive, and $[-\inf_{r\in\R}y(r),\lb_+(c))\to C(\R,\R),\,\lb\mapsto l_\lb$ is continuous in the uniform topology.
\end{teor}
\begin{proof}
As seen before Lemma \ref{lema:order}, \eqref{eq:cond_mu+} (resp.~\eqref{eq:cond_mu-}) is equivalent to the first (resp.~second) inequality in \eqref{eq:desmu}. So, Theorem~\ref{teor:bistability_comparison_autonomous} proves the first two assertions. Now, since $\mu_+(c)<\mu_-(c)$, we deduce from Proposition~\ref{prop:encaje} and \eqref{eq:desmu} that $\lb_-(c)<\mu_+(c)<\mu_-(c)<\lb_+(c)$ if \eqref{eq:cond_mu+} and \eqref{eq:cond_mu-} hold. In this situation, Theorem 4.4 of Ref.~\onlinecite{dno6} proves the remaining assertions of the theorem, excepting two of them. The first one is the existence of exactly three hyperbolic solutions for $\lb\in[\mu_+(c),\mu_-(c)]$, which is ensured by \eqref{eq:desmu} and Theorem~\ref{teor:bistability_comparison_autonomous}. The second one is that, if $\sup_{r\in\R}y(r)-\inf_{r\in\R}y(r)\le\lb_1(c)$, then $\lb_-(c)\ge-\inf_{r\in\R}y(r)$, and hence $\mI_2(c)=(\lb_-(c),\,\lb_+(c))$: in this case, $-\inf_{r\in\R}y(r)\le\lb_1(c)-\sup_{r\in\R}y(r)\le\lb_-(c)$ (see Proposition~\ref{prop:encaje}).
\end{proof}
\begin{coro}\label{cor:cotamala}
If $c>4$,
\begin{equation}\label{eq:cota_no_optima}
\sup_{r\in\R}y(r)-\inf_{r\in\R}y(r)<h_3(c)\,,
\end{equation}
and $\mu_-(c)\ge-\inf_{r\in\R}y(r)$, then all the conclusions of Theorem \ref{teor:concave_lineal_bounds} hold. In addition, $\mI_2(c)=(\lb_-(c),\,\lb_+(c))$ at least for  $c\in(4,\,456]$.
\end{coro}
\begin{proof}
It is clear that the left-hand side of \eqref{eq:cota_no_optima} is greater than or equal to those of \eqref{eq:cond_mu+} and \eqref{eq:cond_mu-}, and Lemma \ref{lema:order} ensures that $h_3(c)<h_2(c)$ for $c>4$.
It remains to check that $\sup_{r\in\R}y(r)-\inf_{r\in\R}y(r)\le\lb_1(c)$ for $c\in(4,\,456]$, which follows from $h_3(c)\le\lb_1(c)$. This last inequality can be proved by squaring twice, which takes the inequality to $p(c)\le0$ for $p(c):=c^4 - 456\,c^3 - 24\,c^2 - 1504\,c + 336$, and then checking that $p$ is non-positive on $(4,\,456]$.
\end{proof}
Observe finally that condition \eqref{eq:cota_no_optima} is more restrictive than \eqref{cond:h1} (see Lemma \ref{lema:order}), but the interval $\mI_2(c)$ of bistability that it provides, that is optimal under the assumed conditions, is larger than the interval $\mI_1(c)$ of Theorem \ref{teor:bistability_comparison_autonomous} (see Proposition \ref{prop:encaje}).
\begin{remark}\label{rm:panorama}
To offer a complete panorama of the situation, we add some of the information obtained in Theorem 4.4 of Ref.~\onlinecite{dno6} under conditions $\lb_-(c)<\mu_+(c)<\mu_-(c)<\lb_+(c)$ in our transitive case. In this nonautonomous situation, it is possible (not sure) the existence of more than one bounded or even hyperbolic attractive solution for $\lb>\lb_+(c)$, but all the solutions asymptotically approach $u_\lb$ as $t\to\infty$. And the same happens with $l_\lb$ to the left of $\lb_-(c)$ if $-\inf_{r\in\R}y(r)<\lb_-(c)$.

In particular, this explains the previous comment about the existence of at least one nonautonomous saddle-node bifurcation point: as $\lb$ approaches $\lb_+(c)$ from the left, the two lower hyperbolic solutions approach each other, giving rise to two non-uniformly-separated and non-hyperbolic solutions (which may coincide) and to the absence of bounded solutions ``close to them" for $\lb>\lb_+(c)$ (see Remark 3.7 in Ref.~\onlinecite{dno6}). The situation is similar for $\tilde\lb_-(c)$ if $\tilde\lb_-(c)>-\inf_{r\in\R}y(r)$, which, as basically proved in Theorem \ref{teor:concave_lineal_bounds}, is the case if $\sup_{r\in\R}y(r)-\inf_{r\in\R}y(r)<\lb_1(c)$. The interested reader can find in the diagrams explaining Theorem 4.5 of Ref.~\onlinecite{dno5} (Figure 6) a depiction of this behavior.
\end{remark}

\section{Precise estimates for periodic and almostperiodic input functions}
\label{sec:precise_estimates}
In this section, we observe that the bounds required in Theorem \ref{teor:concave_lineal_bounds} can be significantly better than those of Corollary \ref{cor:cotamala} in some examples.

\textbf{A trigonometric example.}
We take $y(t)$ as a simple trigonometric function of frequency $\theta$ and amplitude $a>0$, which allows us to get simple expressions for the left-hand terms of \eqref{eq:cond_mu+} and \eqref{eq:cond_mu-}.
Of course, this is not always the case: in general, the condition in Corollary \ref{cor:cotamala} is much easier to check.
\begin{prop}\label{prop:cotas_caso_periodico}
Let $a,\,a_0,\,\phi\in\R$ and $\theta>0$ and
\[
y(t):=a_0+a\,\cos(\theta\,t+\phi)\,.
\]
If $\md>0$, then
\[
\begin{split}
\sup_{r\in\R}\,\md\int_0^\infty e^{-\md\,s}y(r+s)\,ds-\inf_{r\in\R}y(r)&=|a|\left(1+\frac{\md}{\sqrt{\md^2+\theta^2}}\right)\,,\\
\sup_{r\in\R}y(r)-\inf_{r\in\R}\,\md\int_0^\infty e^{-\md\,s}y(r+s)\,ds&=|a|\left(1+\frac{\md}{\sqrt{\md^2+\theta^2}}\right)\,.
\end{split}
\]
\end{prop}
\begin{figure}
     \centering
         \includegraphics[width=\columnwidth]{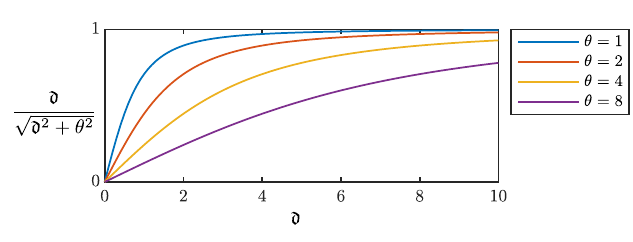}
         \caption{Representation of $\md\mapsto\md/\sqrt{\md^2+\theta^2}$ for different values of $\theta$.}
        \label{fig:periodicexample}
\end{figure}
\begin{proof}
Since the left-hand terms of the equalities of the statement do not depend on $a_0$, we assume without restriction that $a_0=0$.
Let $\mathcal{L}\{y\}(\md):=\int_0^\infty e^{-\md\,s}\,y(s)\,ds$ be the Laplace transform of $y$. Recall that $y{\cdot}r(s):=y(r+s)$.
Since
\[
\begin{split}
\int_0^\infty\!\!\!e^{-\md\,s}\cos(\theta\,(s\!+\!r)\!+\!\phi)\,ds
&=\mathrm{Re}\!\left(e^{i(\theta r+\phi)}\!\int_0^\infty\!\!e^{s(-\md+i\theta)}\,ds\right)
\\[1ex]
&=\frac{\md\,\cos(\theta\,r\!+\!\phi)\!-\!\theta\sin(\theta\,r\!+\!\phi)}{\md^2+\theta^2}\:,
\end{split}
\]
we obtain
\[
\md\,\mathcal{L}\{y{\cdot}r\}(\md)=a\:\frac{\md^2\,\cos(\theta\,r+\phi)-\theta\,\md\,\sin(\theta\,r+\phi)}{\md^2+\theta^2}\,.
\]
It is easy to check that the local extremes of $r\mapsto \md\,\mathcal{L}\{y{\cdot}r\}(\md)$ are attained at the points $r$ with $\tan(\theta\,r+\phi)=-\theta/\md$. This leads to two possibilities: either $\cos(\theta r+\phi)=\md/\sqrt{\md^2+\theta^2}$
and $\sin(\theta r+\phi)=-\theta/\sqrt{\md^2+\theta^2}$, or $\cos(\theta r+\phi)=-\md/\sqrt{\md^2+\theta^2}$
and $\sin(\theta r+\phi)=\theta/\sqrt{\md^2+\theta^2}$.
We evaluate at these points to get
\begin{equation}
\label{eq:max_min_laplace_transform}
\begin{split}
\sup_{r\in[0,2\pi]}\md\,\mathcal{L}\{y{\cdot}r\}(\md)&=|a|\,\frac{\md}{\sqrt{\md^2+\theta^2}}\:,\\
\inf_{r\in[0,2\pi]}\md\,\mathcal{L}\{y{\cdot}r\}(\md)&=-|a|\,\frac{\md}{\sqrt{\md^2+\theta^2}}\:.
\end{split}
\end{equation}
The equalities in the statement follow from here and from $\inf_{r\in\R}y(r)=-|a|$ and $\sup_{r\in\R}y(r)=|a|$.
\end{proof}
Thus, condition \eqref{eq:cond_mu+} is in this case
\begin{equation}
\label{eq:number_asked_1}
|a|\left(1+\frac{\md}{\sqrt{\md^2+\theta^2}}\right)<h_2\big(4(\md+1)\big)\,,
\end{equation}
and condition \eqref{eq:cond_mu-} is
\begin{equation}
\label{eq:number_asked_2}
|a|\left(1+\frac{\md}{\sqrt{\md^2+\theta^2}}\right)<h_3\big(4(\md+1)\big)\,.
\end{equation}
When $\md>0$ is small (i.e., when $c$ is close to $4$), the left-hand sides of these two inequalities differ significantly from $2|a|=\sup_{r\in\R}y(r)-\inf_{r\in\R}y(r)$, which is the expression involved in \eqref{eq:cota_no_optima}.
More precisely, for any given frequency $\theta$, there will be a range of values of $\md>0$ sufficiently small for
which \eqref{eq:cond_mu+} and \eqref{eq:cond_mu-} are significantly less restrictive than \eqref{eq:cota_no_optima}:
see Figure \ref{fig:periodicexample}.
Conversely, for any given frequency $\theta$, as $\md$ becomes sufficiently large, the right-hand side of these bounds approaches $2|a|$.

On the other hand, for a fixed value of $\md>0$, we can always take a sufficiently large frequency $\theta$ so that (31) and (32) are significantly less restrictive than (32). Note that trigonometric terms with high frequencies arise naturally in a Fourier series, as that  giving rise to the almost-periodic example described below.

We emphasize that $\sup_{r\in\R}\md\mathcal{L}\{y{\cdot}r\}(\md)-\inf_{r\in\R}y(r)$ and $\sup_{r\in\R}y(r)-\inf_{r\in\R}\md\mathcal{L}\{y{\cdot}r\}(\md)$ are not, in general, equal. For instance, for $y(t):=(8/25)(2 + \cos(t) - \cos(2t))$, which oscillates from $0$ to $1$, and for $\md=1$,
\[
\begin{split}
\sup_{r\in\R}\mathcal{L}\{y{\cdot}r\}(1)-\inf_{r\in\R}y(r)&\approx0.873\,,\\
\sup_{r\in\R}y(r)-\inf_{r\in\R}\mathcal{L}\{y{\cdot}r\}(1)&\approx0.725\,,
\end{split}
\]
where the maximization and the minimization of the Laplace transform have been carried out numerically.\medskip

\textbf{An almost periodic example.}
In this section, we will use the previous result to establish a bound for the quantities $\sup_{r\in\R}\md\int_0^\infty e^{-\md\,s}y(r+s)\,ds-\inf_{r\in\R}y(r)$ and $\sup_{y\in\R}y(r)-\inf_{r\in\R}\md\int_0^\infty e^{-\md\,s}y(r+s)\,ds$ in the case where the function \( y \) is given by a trigonometric series.
Furthermore, we will demonstrate that this bound is optimal when the frequencies of the trigonometric series are rationally independent, and we will provide a numerical example in which this bound is tighter than $\sup_{r\in\R}y(r)-\inf_{r\in\R}y(r)$.
\begin{prop}
\label{prop:cotas_caso_almostperiodico}
Let $a_0,a_i,\phi_n\in\R$ and $\theta_n>0$ for $i,n\in\N$, let
\begin{equation}
\label{eq:number_asked_3}
y(t):=a_0+\sum_{n=1}^\infty
a_n\,\cos(\theta_n\,t+\phi_n)\,,
\end{equation}
and assume that $\sum_{n=1}^\infty|a_n|$ converges.
If $\md>0$, then
\[
\begin{split}
\sup_{r\in\R}\,\md\!\int_0^\infty\!\!\!e^{-\md\,s}y(r+s)\,ds\!-\!\inf_{r\in\R}y(r)&\!\leq\!\!\sum_{n=1}^\infty\!|a_n|\!\left(\!1\!+\!\frac{\md}{\sqrt{\md^2+\theta_n^2}}\!\right)\,,\\
\sup_{r\in\R}y(r)\!-\!\inf_{r\in\R}\,\md\!\int_0^\infty\!\!\!e^{-\md\,s}y(r+s)\,ds\!&\leq\!\!\sum_{n=1}^\infty\!|a_n|\!\left(\!1\!+\!\frac{\md}{\sqrt{\md^2+\theta_n^2}}\!\right)\,.\\
\end{split}
\]
If, in addition, the set $\{\theta_n\mid\,n\in\N\}$ is rationally independent, then the equalities hold.
Moreover, in this case,
\begin{equation}
\label{eq:sup-inf_quasiperiodic_irrational}
\sup_{r\in\R}y(r)-\inf_{r\in\R}y(r)=2\,\sum_{n=1}^\infty\,|a_n|\,.
\end{equation}
\end{prop}
\begin{proof}
We can assume again $a_0=0$, without restriction. To simplify the notation, we define
\[
y_n(t):=a_n\cos\big(\theta_n\,t+\phi_n\big)\quad\text{and}\quad s_N(t):=\sum_{n=1}^N y_n(t)\,,
\]
for $n\in\N$. Then,
\[
\begin{split}
\sup_{r\in\R}\md\mathcal{L}\{y{\cdot}r\}(\md)&\leq \sum_{n=1}^\infty
\sup_{r\in\R}\,a_n\,\md\mathcal{L}\{y_n{\cdot}r\}(\md)\\
&\leq \sum_{n=1}^\infty\frac{|a_n|\,\md}{\sqrt{\md^2+\theta_n^2}}\:,
\end{split}
\]
where the second inequality is obtained by taking $a=1$
in \eqref{eq:max_min_laplace_transform}.
The first bound of the statement follows from this and
\begin{equation}\label{eq:equality_also_hold}
\inf_{r\in\R}y(r)\ge -\sum_{n=1}^\infty|a_n|\,.
\end{equation}
The second one is analogous.

Let us now prove that the first equality holds under the additional hypothesis that $\{\theta_n\mid\, n\in\N\}$ is rationally independent.
First, let us check that
\begin{equation}\label{eq:proving2}
\sup_{r\in\R}\,\mathcal{L}\{s_N{\cdot}r\}(\md)=\sum_{n=1}^N\,\sup_{r\in\R}\,\mathcal{L}\{y_n{\cdot}r\}(\md)
\end{equation}
for any $N\in\N$.
The inequality $\leq$ is immediate.
Let the sequence $(r_n)$ be such that $\cos(\theta_n r_n+\phi_n)=\md/\sqrt{\md^2+\theta_n^2}$
and $\sin(\theta_n r_n+\phi_n)=-\theta_n/\sqrt{\md^2+\theta_n^2}$.
Then, the proof of Proposition \ref{prop:cotas_caso_periodico} guarantees that
\begin{equation}\label{eq:supremum_at_rn}
\sup_{r\in\R}\mathcal{L}\{y_n{\cdot}r\}(\md)=\mathcal{L}\{y_n{\cdot}r_n\}(\md)\,.
\end{equation}
To prove the inequality $\geq$ in \eqref{eq:proving2}, it suffices to observe that, since $\{\theta_n\mid n\in\N\}$ is rationally independent, the orbit $t\mapsto(\theta_1\,t,\theta_2\,t,\dots,\theta_N\,t)$ is dense in $\T^N$.
So, there exists a sequence $(t_m)$ such that, if $n\in\{1,2,\cdots,\N\}$ then $\lim_{m\to\infty}\theta_n\,t_m=\theta_n\,r_n$ $(\mathrm{mod}\,2\pi)$ and hence $\lim_{m\to\infty}e^{-\md\,s}y_n(s+t_m)=  e^{-\md\,s}y_n(s+r_n)$ uniformly for $s\in[0,\infty)$. This provides
\[
\lim_{n\to\infty}\mathcal{L}\{s_N{\cdot}t_m\}(\md)=\!\sum_{n=1}^N\!\mathcal{L}\{y_n{\cdot}r_n\}(\md)=\!\sum_{n=1}^N\sup_{r\in\R}\mathcal{L}\{y_n{\cdot}r\}(\md)\,,
\]
where the second equality comes from \eqref{eq:supremum_at_rn}, thereby establishing the desired equality in \eqref{eq:proving2}.

Now, let us prove that
\begin{equation}\label{eq:proving3}
\lim_{N\to\infty}\,\sup_{r\in\R}\,\md\,\mathcal{L}\{s_N{\cdot}r\}(\md)=\sup_{r\in\R}\,\md\,\mathcal{L}\{y{\cdot}r\}(\md)\,.
\end{equation}
We fix $\ep>0$ and take $N_0\in\N$ such that $\sum_{n=N_0+1}^\infty|a_n|<\ep$.
Then,
\[
\begin{split}
&\left|\,\sup_{r\in\R}\,\md\,\mathcal{L}\{y{\cdot}r\}(\md)-\sup_{r\in\R}\,\md\,\mathcal{L}\{s_{N_0}{\cdot}r\}(\md)\,\right|\\
&\qquad\leq\left|\,\sup_{r\in\R}\int_0^\infty \md\,e^{-\md\,s}\!\!\sum_{n=N_0+1}^\infty a_n\cos\big(\theta_n(s+r)+\phi_n\big)\,\right|<\ep\,,
\end{split}
\]
which ensures \eqref{eq:proving3}.
It then suffices to take the limit as $N\to\infty$ in \eqref{eq:proving2} and apply \eqref{eq:proving3} and \eqref{eq:max_min_laplace_transform} to obtain
\begin{equation}\label{eq:la_anterior_asked}
\sup_{r\in\R}\,\mathcal{L}\{y{\cdot}r\}(\md)=\sum_{n=1}^\infty\frac{|a_n|\,\md}{\sqrt{\md^2+\theta_n^2}}\,.
\end{equation}
An analogous but simpler $\T^N$-density argument shows that the equality also holds in \eqref{eq:equality_also_hold}.
Combining this and \eqref{eq:la_anterior_asked} provides the equality in the statement.

The arguments to obtain the equality on the second inequality of the statement and on \eqref{eq:sup-inf_quasiperiodic_irrational} in this case are analogous.
\end{proof}
\begin{figure}
     \centering
         \includegraphics[width=\columnwidth]{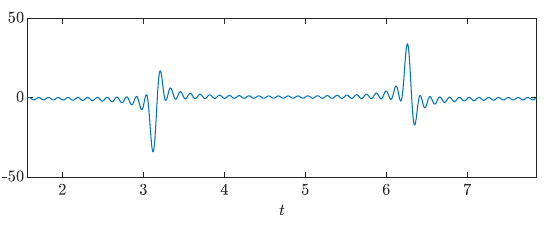}
         \caption{Numerical depiction of the $2\pi$-periodic function $y(t)=\sum_{n=1}^{51}(\cos(n\,t)-\sin(n\,t))$ on $[\pi/2,5\pi/2]$.}
        \label{fig:ejemploalmostperiodic1}
\end{figure}
\begin{figure}
     \centering
         \includegraphics[width=\columnwidth]{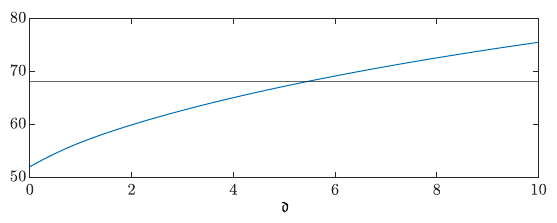}
         \caption{Depiction of the curve $\md\mapsto\sum_{n=1}^{51}2(1+\md/\sqrt{\md^2+n^2})$, which is the value of the right-hand side of the bounds in Proposition \ref{prop:cotas_caso_almostperiodico} for $y(t)=\sum_{n=1}^{51}(\cos(n\,t)-\sin(n\,t))$.
         The horizontal line represents $68.03$, which is an approximation of $\sup_{r\in\R}y(r)-\inf_{r\in\R}y(r)$.}
        \label{fig:ejemploalmostperiodic2}
\end{figure}
The most widespread example in which $\{\theta_n\mid\,n\in\N\}$ is not rationally independent is that of Fourier series, which we will use in what follows.
Let us show with a particular example of the map \eqref{eq:number_asked_3} that the bound provided by the right-hand term of the inequalities in Proposition \ref{prop:cotas_caso_almostperiodico} can be below $\sup_{r\in\R}y(r)-\inf_{r\in\R}y(r)$ for small values of $\md>0$.
So, changing the left-hand term in \eqref{eq:number_asked_1} and \eqref{eq:number_asked_2} by the right-hand term of the inequalities in Proposition \ref{prop:cotas_caso_almostperiodico}, we obtain easily computable conditions to ensure bistability if $c>4$ is not too large, whereas for larges values of $c$, the bound provided by $\sup_{r\in\R}y(r)-\inf_{r\in\R}y(r)$ is better. We take
\[
y(t):=\sum_{n=1}^{51}\big(\cos(n\,t)-\sin(n\,t)\big)\,,
\]
depicted in Fig.~\ref{fig:ejemploalmostperiodic1}, for which $\sup_{r\in\R}y(r)-\inf_{r\in\R}y(r)\approx 68.03$.
As shown in Fig.~\ref{fig:ejemploalmostperiodic2}, for values of $d$ below a certain threshold, the bound proposed in Proposition \ref{prop:cotas_caso_almostperiodico} is lower than $68.03$.

Proposition \ref{prop:cotas_caso_almostperiodico} requires the absolute convergence of the series of coefficients of the function that captures the temporal variation of the input. However, this is not always the case, not even for the Fourier series of a periodic function. The following proposition employs Ces\`{a}ro sums to establish a sufficient condition less restrictive than their absolute convergence.
\begin{prop} Let $y(t)$ be a continuous $2\pi$-periodic function with Fourier series
\begin{equation}\label{eq:fourierseries}
a_0+\sum_{n=1}^\infty\big(a_n\cos(nt)+b_n\sin(nt)\big)\,.
\end{equation}
Then, the two inequalities of Proposition {\rm \ref{prop:cotas_caso_almostperiodico}} hold when replacing the right-hand side with
\[
\lim_{N\to\infty}\,\frac{1}{N}\sum_{n=1}^{N-1}(N-n)\big(|a_n|+|b_n|\big)\left(1+\frac{\md}{\sqrt{\md^2+n^2}}\right).
\]
\end{prop}
\begin{proof}
Once again, we can assume $a_0=0$. Let $s_n$ be the $n$-partial sum of the Fourier series \eqref{eq:fourierseries} for $n\ge 1$.
The $N$-th Ces\`{a}ro mean of this series is
 \[
\begin{split}
\sigma_N(t)&:=\frac{1}{N}\sum_{n=0}^{N-1}s_n(t)\\
&=\frac{1}{N}\sum_{n=1}^{N-1}(N-n)\big(a_n\cos(nt)+b_n\sin(nt)\big)
\end{split}
\]
for $N\ge 2$.
Theorem 3.1 of Ref.~\onlinecite{katznelson1978}, which is due to Fej\'{e}r, ensures that $\lim_{N\to\infty}\sigma_N(t)=y(t)$ uniformly on $\R$.

On the other hand,
\[
\begin{split}
&\sup_{r\in\R}\md\int_0^\infty e^{-\md\,s}\sigma_N(r+s)\,ds\\
&\quad\le\frac{\md}{N}\!\sum_{n=1}^{N-1}
(N-n)
\left(
\sup_{r\in\R}\,a_n\int_0^\infty\! e^{-\md\,s}\cos\big(n(r+s)\big)\,ds\right.\\
&\left.\quad\qquad\qquad\qquad\qquad+\sup_{r\in\R}\,b_n\int_0^\infty\! e^{-\md\,s}\sin\big(n(r+s)\big)\,ds\right),
\end{split}
\]
so the calculations in the proof of Proposition \ref{prop:cotas_caso_periodico} ensure that
\begin{equation}\label{eq:inprogress}
\begin{split}
&\sup_{r\in\R}\md\int_0^\infty e^{-\md\,s}\sigma_N(r+s)\,ds-\inf_{r\in\R}\sigma_N(r)\\
&\qquad\leq\frac{1}{N}\sum_{n=1}^{N-1}(N-n)\big(|a_n|+|b_n|\big)\left(1+\frac{\md}{\sqrt{\md^2+n^2}}\right),
\end{split}
\end{equation}
since
\[
-\inf_{r\in\R}\sigma_N(r)\leq \frac{1}{N}\sum_{n=1}^{N-1}(N-n)\big(|a_n|+|b_n|\big)\,.
\]
An analogous calculation can be carried out with
\begin{equation}\label{eq:number_asked_4}
\sup_{r\in\R}\sigma_N(r)-\inf_{r\in\R}\,\md\!\int_0^\infty e^{-\md\,s}\sigma_N(r+s)\,ds\,.
\end{equation}
The uniform convergence $\sigma_N(t)\to y(t)$ as $N\to\infty$ ensures that
\[
\lim_{N\to\infty}\,\md\!\int_0^\infty e^{-\md\,s}\sigma_N(r+s)\,ds=\md\!\int_0^\infty e^{-\md\,s}y(r+s)\,ds
\]
uniformly for $r\in\R$.
Therefore,
\[
\lim_{N\to\infty}\sup_{r\in\R}\,\md\!\int_0^\infty e^{-\md\,s}\sigma_N(r+s)\,ds= \sup_{r\in\R}\,\md\!\int_0^\infty e^{-\md\,s}y(r+s)\,ds\,.
\]
We can reason in the same way with $\inf_{r\in\R}\sigma_N(r)$ to get the first inequality indicated in the statement, and proceed similarly with the inequality obtained for \eqref{eq:number_asked_4}.
\end{proof}
\section{Further results on bistability: the auxiliary d-concavity band}
\label{sec:bistability_dconcavity}
This section provides additional insight into the occurrence of bistability by leveraging the existence of the d-concavity band surrounding the graph of the inflection curve $x=\sqrt3$, which is determined by
the interval $\big[\sqrt{3-2\sqrt{2}},\sqrt{3+2\sqrt{2}}\big]$ (see Section \ref{sec:bonifacio-lugiato}).
For $c\in(4,\,2+2\sqrt{2})\approx(4,\,4.82842)$, we describe a scenario similar to that of Corollary \ref{cor:cotamala} under different conditions on the variation of $y$ which are significantly less restrictive for $c$ close to $4$. The next new values of the parameter, determined by the endpoints of the d-concavity band, play a role on the result:
\[
\begin{split}
\lb_3(c)&:=\!-\bar g\!\left(\!\sqrt{3-2\sqrt{2}}\right)\!=\!\sqrt{3-2\sqrt{2}}\left(\!\!\left(1+\frac{\sqrt{2}}{2}\right)c+1\right),\\
\lb_4(c)&:=\!-\bar g\!\left(\!\sqrt{3+2\sqrt{2}}\right)\!=\!\sqrt{3+2\sqrt{2}}\left(\!\!\left(1-\frac{\sqrt{2}}{2}\right)c+1\right).
\end{split}
\]
Note that $c\mapsto\lb_3(c)$ and $c\mapsto\lb_4(c)$ are globally defined linear maps with the same slope, and hence $\lb_4(c)-\lb_3(c)=\lb_4(0)-\lb_3(0)=2$.
Recall that $\lb_1(c)$ and $\lb_2(c)$ are defined by \eqref{eq:definition_lb1_and_lb2}.
\begin{teor}\label{teor:dconcavity}
Let $c\in\mC:=(4,\,2+2\sqrt{2})$ be fixed. If $\sup_{r\in\R}y(r)-\inf_{r\in\R}y(r)<2=\lb_4(c)-\lb_3(c)$ and
\begin{equation}\label{eq:intervalo}
 \lb\in\left(\lb_3(c)-\inf_{r\in\R}y(r)\,,\;\lb_4(c)-\sup_{r\in\R}y(r)\right),
\end{equation}
then all the bounded positive solutions of \eqref{eq:parametricnonautonomousbonifacio}$_\lb$ take values in the d-concavity band of $g$, $\big[\sqrt{3-2\sqrt{2}},\sqrt{3+2\sqrt{2}}\,\big]$.
In addition, $h_4(c):=\lb_4(c)-\lb_1(c)$ is strictly positive and, if
\begin{equation}\label{eq:min}
 \sup_{r\in\R}y(r)-\inf_{r\in\R}y(r)<\min(h_1(c),\,h_4(c))\,,
\end{equation}
then there exists an interval $\mI_3(c):=(\tilde\lb_1(c),\,\tilde\lb_2(c))\cap[-\inf_{r\in\R}y(r),\,\tilde\lb_2(c))$ such that \eqref{eq:parametricnonautonomousbonifacio}$_\lb$ exhibits bistability
for $\lb\in\mI_3(c)$ and not for any other $\lb\geq-\inf_{r\in\R}y(r)$.

More precisely: for $\lb\in\mI_3(c)$ there exist three positive hyperbolic solutions of \eqref{eq:parametricnonautonomousbonifacio}$_\lb$, $l_\lb<m_\lb<u_\lb$, with $l_\lb$ and $u_\lb$ attractive and with basins of attraction on the positive half-plane separated by the graph of $m_\lb$, which is repulsive;
if $\lb\ge\tilde\lb_2(c)$, then $u_\lb$ is the unique positive hyperbolic solution, it is attractive, and $(\tilde\lb_1(c),\,\infty)\cap[-\inf_{r\in\R}y(r),\,\infty)\to C(\R,\R),\,\lb\mapsto u_\lb$ is continuous in the uniform topology; and if $\tilde\lb_1(c)\ge-\inf_{r\in\R}y(r)$ and $\lb\in[-\inf_{r\in\R}y(r),\,\tilde\lb_1(c)]$, then $l_\lb$ is the unique positive hyperbolic solution, it is attractive, and $[-\inf_{r\in\R}y(r),\,\tilde\lb_2(c))\to C(\R,\R),\,\lb\mapsto l_\lb$ is continuous in the uniform topology. Moreover, for $\lb>\tilde\lb_2(c)$ and for $\lb\in[-\inf_{r\in\R}y(r),\,\tilde\lb_1(c))$ (if nonempty), every solution approaches the unique attractive hyperbolic solution as time increases.
\end{teor}
\begin{proof}
It is easy to deduce from $c\in\mC$ that
\begin{equation}\label{eq:contencion_cinmC}
 [\,x_1(c),\:x_2(c)\,]\subset\left( \!\sqrt{3-2\sqrt{2}},\,\sqrt{3+2\sqrt{2}}\,\right)
\end{equation}
(see \eqref{eq:definition_x1_and_x2}). In fact, the value $c_0:=2+2\,\sqrt{2}$ is optimal to get $x_2(c)<\sqrt{3+2\sqrt{2}}$ for $c\in\mC$: recall that $x_2$ is strictly increasing on $[4,\infty)$, and observe that $x_2(c_0)=\sqrt{3+2\,\sqrt{2}}$.

We fix $\lb$ satisfying \eqref{eq:intervalo}. Recall that $\bar g$ is strictly decreasing outside $[\,x_1(c),\,x_2(c)]$ (see the beginning of Section \ref{sec:bistability_autonomous}). This ensures, first, that the invariance condition $\lb\ge-\inf_{r\in\R}y(r)$ holds, since $\lb+\inf_{r\in\R}y(r)>\lb_3(c)>-\bar g(0)=0$; second, that $\lb+y(t)+\bar g(r)>\lb+\inf_{r\in\R}y(r)-\lb_3(c)=:\delta^1_\lb>0$ if $r<\sqrt{3-2\sqrt{2}}$, which precludes the existence of solutions of $x'=\lb+y(t)+\bar g(x)$ with bounded backward semiorbit taking any value on $\big(\!-\infty,\sqrt{3-2\sqrt{2}}\big)$; and third, that $\lb+y(t)+\bar g(r)<\lb+\sup_{r\in\R}y(r)-\lb_4(c)=:-\delta^2_\lb<0$ if $r>\sqrt{3+2\sqrt{2}}$, which precludes the existence of solutions with bounded backward semiorbit of $x'=\lb+y(t)+\bar g(x)$ taking any value on $\big(\!\sqrt{3+2\sqrt{2}},\infty\big)$. So, the assertion about the bounded positive solutions of \eqref{eq:parametricnonautonomousbonifacio}$_\lb$ follows from $g=\bar g$ on $x\ge 0$ and from the invariance of the positive half-plane for \eqref{eq:parametricnonautonomousbonifacio}$_\lb$.

Condition \eqref{eq:min} involves two requirements. The first one, $\sup_{r\in\R}y(r)-\inf_{r\in\R}y(r)<h_1(c)$ has two consequences. First, the non-emptyness of the interval $\mI_1(c)$ of Theorem \ref{teor:bistability_comparison_autonomous} of values of $\lb$ for which \eqref{eq:parametricnonautonomousbonifacio}$_\lb$ has three hyperbolic solutions. Second, the non-emptyness of the interval of \eqref{eq:intervalo}, which follows from $h_1(c)<2$ for $c\in\mC$. This property, which is clearly reflected in Figure \ref{fig:j_functions}, can be proved by squaring twice, which takes the inequality to $p(c):=c^4-12\,c^3+40\,c^2-96\,c+32<0$. It is an elementary exercise to prove that the polynomial $p$ has exactly two zeros, one in $(0,1)$ and the other in $(9,10)$, and that it is negative in $\mC$. The second requirement in \eqref{eq:min}, $\sup_{r\in\R}y(r)-\inf_{r\in\R}y(r)<h_4(c)$, is intended to guarantee that the intervals of \eqref{eq:intervalo} and $\mI_1(c)$ have nonempty intersection, what we prove in what follows. Since $x_2(c)<\sqrt{3+2\sqrt{2}}$ and $x_2(c)$ is the maximum of $\bar g$ on $[x_2(c),\infty)$, $\lb_4(c)=-\bar g\big(\sqrt{3+2\sqrt{2}}\big)>-\bar g(x_2(c))=\lb_1(c)$ for all $c\in\mC$. On the other hand,
$\lb_3(c)<\lb_1(c)$, as deduced from the increasing character of both maps: $\lb_3(c)<\lb_3(2+2\,\sqrt2)<4<\lb_1(4)\le\lb_1(c)$ for all $c\in\mC$. So, we get $\lb_3(c)-\inf_{r\in\R}y(r)<\lb_1(c)-\inf_{r\in\R}y(r)<\lb_4(c)-\sup_{r\in\R}y(r)$ for $c\in\mC$ if $\sup_{r\in\R}y(r)-\inf_{r\in\R}y(r)<\lb_4(c)-\lb_1(c)=h_4(c)$, and this proves the assertion.

Hence, if \eqref{eq:min} holds, and if $\lb_0$ satisfies both \eqref{eq:intervalo} and \eqref{eq:lambda_interval_autonomous_prop},
then Theorem~\ref{teor:bistability_comparison_autonomous} and the property guaranteed by \eqref{eq:intervalo} ensure the existence of three hyperbolic solutions $l_{\lb_0},\,m_{\lb_0}$ and $u_{\lb_0}$ of \eqref{eq:concave-linear_bonifacio_bifurcation}$_{\lb_0}$ with
$\sqrt{3-2\sqrt{2}}\leq l_{\lb_0}(t)<x_1(c)<m_{\lb_0}(t)<x_2(c)<u_{\lb_0}(t)\leq\sqrt{3+2\sqrt{2}}$ for all $t\in\R$.
\par
Let $\tilde g$ be a globally d-concave $C^2$ extension of $\bar g$ coinciding with it at $\big[\sqrt{3-2\sqrt{2}},\,\sqrt{3+2\sqrt{2}}\big]$ and satisfying $\lim_{x\to\pm\infty}\tilde g(x)=-\infty$, which can be constructed by taking third degree polynomials outside the interval: see, e.g., the proof of Theorem 4.11 in Ref.~\onlinecite{dno6}. It is not hard to check that the family $x'=\lb+\w(t)+\tilde g(x)$, $\w\in\W_y$, satisfies all the conditions required at Theorem 4.5 of Ref.~\onlinecite{dno5}. Recall also that the $\sigma$-orbit of $y$ is dense in $\W_y$. This result provides an interval $(\tilde\lb_1(c),\tilde\lb_2(c))$ such that:
\begin{itemize}[leftmargin=8pt]
\item[-] for any $\lb>\tilde\lb_1(c)$ (resp.~$\lb<\tilde\lb_2(c)$), the upper (resp.~lower) bounded solution $\tilde u_\lb$ (resp.~$\tilde l_\lb$) of $x'=\lb+y(t)+\tilde g(x)$ is hyperbolic attractive, and the map $(\tilde\lb_1(c),\infty)\to C(\R,\R),\,\lb\mapsto \tilde u_\lb$ (resp.~$(-\infty,\lb_2)\to C(\R,\R),\,\lb\mapsto \tilde l_\lb$) is continuous in the uniform topology and strictly increasing;
\item[-] for $\lb\in(\tilde\lb_1(c),\tilde\lb_2(c))$, the basins of attraction on the positive half-plane of $\tilde l_\lb$ and $\tilde u_\lb$ are separated by the graph of the unique remaining hyperbolic solution, $\tilde m_\lb$, which is repulsive, and the map $(\tilde\lb_1(c),\tilde\lb_2(c))\to C(\R,\R),\,\lb\mapsto \tilde m_\lb$ is continuous in the uniform topology and strictly decreasing;
\item[-] for any $\lb>\tilde\lb_2(c)$ (resp.~$\lb<\tilde\lb_1(c)$), $\lim_{t\to\infty}(\tilde u_\lb(t)-\tilde x_\lb(t))=0$ (resp.~$\lim_{t\to\infty}(\tilde l_\lb(t)-\tilde x_\lb(t))=0$) for any other solution $\tilde x_\lb$ of $x'=\lb+y(t)+\tilde g(x)$;
\item[-] $\tilde l_{\tilde\lb_1(c)}$ is the unique hyperbolic solution of $x'=\tilde\lb_1(c)+y(t)+\tilde g(x)$, it is attractive, and it is uniformly separated from the solutions $\tilde m_{\tilde\lb_1(c)}\le \tilde u_{\tilde\lb_1(c)}$ given by $\tilde m_{\tilde\lb_1(c)}:=\lim_{\lb\to(\tilde\lb_1(c))^+}\tilde m_\lb$ and $\tilde u_{\tilde\lb_1(c)}:=\lim_{\lb\to(\tilde\lb_1(c))^+}\tilde u_\lb$, which are not uniformly separated; and
   $\tilde u_{\tilde\lb_2(c)}$ is the unique hyperbolic solution of $x'=\tilde\lb_2(c)+y(t)+\tilde g(x)$, it is attractive, and it is uniformly separated from the solutions $\tilde m_{\tilde \lb_2}\ge\tilde l_{\tilde\lb_1(c)}$ given by $\tilde m_{\tilde\lb_2(c)}:=\lim_{\lb\to(\tilde\lb_2(c))^-}\tilde m_\lb$ and $\tilde l_{\tilde\lb_2(c)}:=\lim_{\lb\to(\tilde\lb_2(c))^-}\tilde l_\lb$, which are not uniformly separated.
\end{itemize}
In addition, although these properties are not included in the statement of Theorem 4.5 of Ref.~\onlinecite{dno5}, standard comparison arguments can be used to check that
$\lim_{t\to\infty}(\tilde l_{\tilde\lb_1(c)}(t)-\tilde x_{\tilde\lb_1(c)}(t))=0$ if $\tilde x_{\tilde\lb_1(c)}<\tilde m_{\tilde\lb_1(c)}$ and
$\lim_{t\to\infty}(\tilde u_{\tilde\lb_2(c)}(t)-\tilde x_{\tilde\lb_2(c)}(t))=0$ if $\tilde x_{\tilde\lb_2(c)}>\tilde m_{\tilde\lb_2(c)}$.

According to the previous description, the two upper hyperbolic solutions for $\lb\in(\tilde\lb_1(c),\lb_0]$ with $\lb\ge-\inf_{r\in\R}y(r)$ and the two lower hyperbolic solutions for $\lb\in[\lb_0,\tilde\lb_2(c))$ take values in $\big[\sqrt{3-2\sqrt{2}},\,\sqrt{3+2\sqrt{2}}\big]$, they are also hyperbolic solutions for \eqref{eq:parametricnonautonomousbonifacio}$_\lb$.
Let us focus on $\lb\ge\lb_0$. First, we take $\lb\in[\lb_0,\tilde\lb_2(c)]$, call $l_\lb:=\tilde l_\lb$ and $m_\lb:=\tilde m_\lb$, and observe that they solve \eqref{eq:parametricnonautonomousbonifacio}$_\lb$, that they are hyperbolic except for $\lb=\tilde\lb_2(c)$, and that there are no more positive hyperbolic solutions below $m_\lb$.
Let $x_\lb$ denote any solution with $x_\lb>m_\lb$ of \eqref{eq:parametricnonautonomousbonifacio}$_\lb$. Recall that there exists a unique attractive hyperbolic solution $u_\lb$ above $x_2(c)$ attracting any solution taking any value in $[x_2(c),\infty)$ as time increases: see Theorem~ \ref{teor:bistability_comparison_autonomous}. We assume for contradiction that $x_\lb(t)<x_2(c)$ for all $t\in\R$, observe that $x_\lb$ solves $x'=\lb+y(t)+\tilde g(x)$, and get the contradiction from $\lim_{t\to\infty}(\tilde u_\lb(t)-x_\lb(t))=0$ and $\tilde u_\lb>\tilde u_{\lb_0}>x_2(c)$. So, $\lim_{t\to\infty}(u_\lb(t)-x_\lb(t))=0$. This proves the bistability for $\lb\in[\lb_0,\tilde\lb_2(c))$ and the lack of it for $\tilde \lb_2(c)$. Now, we take $\lb>\tilde\lb_2(c)$ and use also the previous description and the information provided by Theorem~\ref{teor:bistability_comparison_autonomous} to conclude that $u_\lb$ is the unique hyperbolic solution as well as the stated asymptotic behavior of the remaining solutions. We also call $u_\lb:=\tilde u_\lb$ for $\lb\in(\tilde \lb_1,\,\tilde\lb_2(c))\cap[-\inf_{r\in\R}y(r),\,\tilde\lb_2(c))$ (which is coherent with the notation of Theorem~\ref{teor:bistability_comparison_autonomous}). The continuity of $(\tilde\lb_1(c),\,\infty)\cap[-\inf_{r\in\R}y(r),\,\infty)\to C(\R,\R),\,\lb\mapsto u_\lb$ in the uniform topology is a consequence of the persistence of hyperbolic solutions (see, e.g., Theorem 3.8 of Ref.~\onlinecite{potzsche2011}). A similar analysis to the left of $\lb_0$---and always to the right of $-\inf_{r\in\R}y(t)$---completes the proof.
\end{proof}
The last result of this section provides bounds for the points $\tilde \lb_1(c)$ and $\tilde\lb_2(c)$ of Theorem \ref{teor:dconcavity} and shows that the interval $\mI_3(c)$ of Theorem \ref{teor:dconcavity} contains the interval $\mI_1(c)$ of Theorem \ref{teor:bistability_comparison_autonomous} under the conditions of the last of the two theorems. Note also that the comment made in Remark \ref{rm:panorama} about the occurrence of one or two saddle-node bifurcation values of hyperbolic solutions also applies to the situation described in Theorem \ref{teor:dconcavity}.
\begin{figure}
     \centering
         \includegraphics[width=\columnwidth]{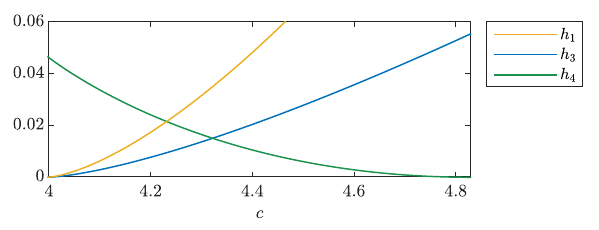}
         \caption{Representation of the functions $h_1(c)$, $h_3(c)$ and $h_4(c)$ defined in \eqref{def:h1}, \eqref{def:h3} and Theorem \ref{teor:dconcavity} on $[4,2+2\sqrt{2}]$.}
        \label{fig:j_functions_2}
\end{figure}
\begin{prop} \label{prop:encajedos}
If $4<c<2+2\sqrt{2}$ and \eqref{eq:min} holds, then
\[
\begin{split}
  \lb_3(c)-\sup_{r\in\R}y(r)&\leq\tilde\lb_1(c)\leq\lb_1(c)-\inf_{r\in\R} y(r)\,,\\
  \lb_2(c)-\sup_{r\in\R}y(r)&\leq\tilde\lb_2(c)\leq\max\big(\lb_2(c),\lb_4(c)\big)-\inf_{r\in\R} y(r)\,.
\end{split}
\]
\end{prop}
\begin{proof}
Again, let $\mC:=(4,\,2+2\sqrt{2})$. The properties---the graph---of $\bar g$ show that it attains its maximum value on $\big[\sqrt{3-2\sqrt{2}},\,\sqrt{3+2\sqrt{2}}\,\big]$ either at $\sqrt{3-2\sqrt{2}}$ or at $x_2(c)$, and we have checked in the previous proof that $\bar g(\sqrt{3-2\sqrt{2}})=-\lb_3(c)>-\lb_1(c)=\bar g(x_2(c))$ for $c\in\mC$. So, $\lb_3(c)+\bar g(x)\le 0$ for $x\in[\sqrt{3-2\sqrt{2}},\sqrt{3+2\sqrt{2}}]$ and $c\in\mC$.
So, if $\lb\le\lb_3(c)-\sup_{r\in\R}y(r)-\delta$ for a $\delta>0$, then any solution of \eqref{eq:parametricnonautonomousbonifacio}$_\lb$ taking all its values in $[\sqrt{3-2\sqrt{2}},\sqrt{3+2\sqrt{2}}]$ would satisfy $x'(t)\le-\delta<0$ for all $t\in\R$, which is impossible; and hence $\tilde\lb_1(c)\geq\lb_3(c)-\sup_{r\in\R}y(r)$, since \eqref{eq:parametricnonautonomousbonifacio}$_{\tilde\lb_1(c)}$ has such solutions. This proves the first inequality for $\tilde\lb_1(c)$. The second inequality for $\tilde\lb_2(c)$ is checked similarly, using $\bar g(x)\ge-\max(\lb_2(c),\lb_4(c))$ for $x\in[\sqrt{3-2\sqrt{2}},\sqrt{3+2\sqrt{2}}]$ and $c\in\mC$.

To check the first inequality for $\tilde\lb_2(c)$,
we combine that $\lb\mapsto l_\lb$ is not right-continuous at $\tilde\lb_2(c)$---proved in Theorem 4.5 of Ref.~\onlinecite{dno5}---with Proposition 4.3 of Ref.~\onlinecite{dno6} and with Proposition~\ref{prop:encaje} to get $\tilde\lb_2(c)\geq\lb_+(c)\ge\lb_2(c)-\sup_{r\in\R}y(r)$. An analogous argument shows that $\tilde\lb_1(c)\leq\lb_-(c)\le\lb_1(c)-\inf_{r\in\R}y(r)$.
\end{proof}
Figure \ref{fig:j_functions_2} shows that, for $c\in(4,4.23)$, conditions \eqref{eq:min} and \eqref{eq:lambda_interval_autonomous_prop} are equivalent, so that for those values of $c$, Theorem \ref{teor:dconcavity} provides additional information beyond that of Theorem \ref{teor:bistability_comparison_autonomous}. It also shows that \eqref{eq:min} is less restrictive than \eqref{eq:cota_no_optima} for an interval of values of $c$ containing $(4,4.32)$. Note finally that $\mI_1(c)\subseteq\mI_2(c)=\mI_3(c)$
if $c\in(4,\,2+2\sqrt{2})$ and \eqref{eq:min} and \eqref{eq:cota_no_optima} are fulfilled.

\section{Relaxation oscillations}\label{sec:relxationoscillations}
We will work with a fixed value of $c>4$. This section analyzes the equation
\begin{equation}\label{eq:yeps}
 x'=y_\ep(t) + \bar g(x)\,,
\end{equation}
obtained by replacing the input $\lb+y(t)$ with a specific map $y_\ep(t)$ depending on two parameters, $\ep>0$ and $r>0$.
The interval image of this map contains $[\lb_1(c),\,\lb_2(c)]$ in its interior (see \eqref{eq:definition_lb1_and_lb2}), so that the results of Section \ref{sec:bistability_autonomous} do not ensure bistability (see Remark \ref{rm:nueva}.2). However, we will find bistability for some pairs $(r,\ep)$, as well as relaxation oscillations---and uniform stability---for other pairs. The arguments to be used are classical within the theory of fast-slow systems (see, e.g., Refs. \onlinecite{berglund2006}, \onlinecite{kuehn2015} and \onlinecite{jones1995}), and are due, among others, to Tikhonov, Pontryagin and Fenichel.\cite{tikhonov1952,pontryagin1957,fenichel1979}

More precisely, we take
\begin{equation}
\label{eq:y(t)_relaxation_oscillations}
y_\ep(t):=\alpha+(\beta+\ep^r)\,\sin(\ep\,t)\,,
\end{equation}
for $\ep>0$ and $r>0$, where
\begin{equation}
\label{eq:definition_alpha_beta}
\alpha:=\frac{\lb_1(c)+\lb_2(c)}{2}\quad\text{and}\quad\beta:=\frac{\lb_2(c)-\lb_1(c)}{2}\,.
\end{equation}
Our analysis will fix $r$ and let $\ep$ vary, which is the reason for the chosen notation. In fact, we are only interested in the dynamical behavior when $\ep>0$ is small, that is, when the time variation of the input $y_\ep$ is slow.

Note that if $\bar y_\ep(t):=\alpha+\beta\sin(\ep\,t)$ were taken instead of $y_\epsilon$, then Theorem~\ref{teor:bistability_comparison_autonomous} would guarantee the existence of bistability for all values of $\epsilon>0$,
since $\lb=0$ would fulfill $\lb\in[\lb_1(c)-\inf_{r\in\R}y(r),\,\lb_2(c)-\sup_{r\in\R}y(r)]$ and there are no constant maps in the hull of $\bar y_\ep$.

The input $y_\ep$ follows a sinusoidal pattern that spans the set $[\lb_1(c)-\ep^r,\,\lb_2(c)+\ep^r]$, which is slightly larger than the closure $[\lb_1(c),\,\lb_2(c)]$ of the set of values of bistability of the autonomous problem \eqref{eq:autonomousbonifacio} (see Fig.~\ref{fig:bif_diagram_aut_bonifacio}). As long as $y_\ep(t)$ remains within the interval of autonomous bistability, there will exist solutions that track their stable equilibria, in accordance with Tikhonov's Theorem.
Consequently, the magnitude and velocity at which the input $y_\ep(t)$ surpasses the autonomous bifurcation points $\lb_1(c)$ and $\lb_2(c)$ will determine whether all the solutions of the equation eventually remain close to one of the stable autonomous branches or they transit from one stable branch to the other. That is, this depends on the sizes of $\ep$ and $r$.

\begin{figure}
     \centering
         \includegraphics[width=0.6\columnwidth]{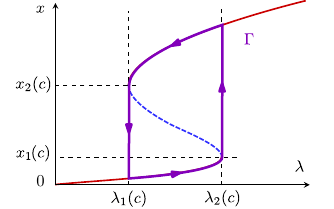}
         \caption{In purple, the limit curve $\Gamma$ of the relaxation oscillations for any $c>4$: in this drawing, $c=5$.}
        \label{fig:bif_diagram_relaxation_definition}
\end{figure}

The autonomous bifurcation diagram depicted in Fig.~\ref{fig:bif_diagram_aut_bonifacio}
leads to the construction of the curve $\Gamma$ to which relaxation solutions---to be defined below---will refer, which is depicted in Fig.~\ref{fig:bif_diagram_relaxation_definition}: it consists of the non-repulsive fixed points of \eqref{eq:autonomousbonifacio} between $\lb_1(c)$ and $\lb_2(c)$ (slow dynamics) together with two vertical lines connecting each of the saddle-node bifurcation points with the other stable branch (fast dynamics).
\begin{figure}
     \centering
         \includegraphics[width=\columnwidth]{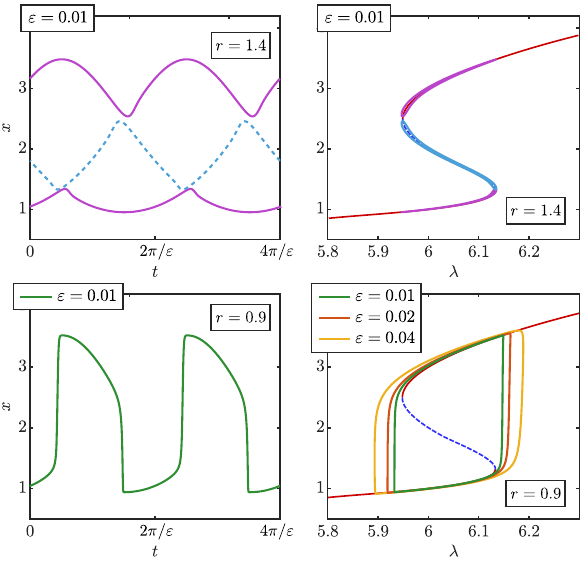}
         \caption{Numerical simulation of the two possible scenarios in \eqref{eq:yeps} with $c=5$ and small $\ep$.
In the top-left panel, the graph of the unique hyperbolic ($2\pi/\ep$-periodic) solution $x_\ep$ for one of these small values of $\ep$ and for $t\in[0,4\pi/\ep]$.
The top-right panel represents the relaxation-oscillation cycles $(y_\ep(t),\, x_\ep(t))$ approaching $\Gamma$, in different colors for different values of $\ep$.
The equilibria of the autonomous system \eqref{eq:autonomousbonifacio} are depicted as in Fig.~\ref{fig:bif_diagram_aut_bonifacio}.
The lower panels illustrate the case of bistability for $r = 1.4$. The representations are analogous to those above, but here for a single value of $\ep$: the attractive solutions are plotted in purple, and the repulsive one in cyan blue.
}
        \label{fig:relaxation_simulation}
\end{figure}

\begin{defi} An $\ep$-family $x_\ep(t)$ of $2\pi/\ep$-periodic solutions of \eqref{eq:yeps} is said to be a \textit{relaxation oscillation} if the closed curves $\{(y_\ep(t),\,x_\ep(t))\mid\,t\in\R\}$ converge in Hausdorff distance to the curve $\Gamma$ as $\ep\to0^+$.
\end{defi}

Let us formulate the last result of this paper.
\begin{teor}
\label{teor:relaxation_oscillations}
Let $c>4$ be fixed.
\begin{itemize}
\item[\rm(i)] If $r>1$, then \eqref{eq:yeps} does not have relaxation oscillations.
In particular, for sufficiently small $\ep>0$, there are exactly three periodic positive solutions: the upper and lower ones are hyperbolic attractive, while the middle one is hyperbolic repulsive and separates the basins of attraction of the other two in the positive half-plane.
That is, \eqref{eq:yeps} exhibits bistability.
\item[\rm(ii)]If $r<1$, then \eqref{eq:yeps} has a relaxation oscillation.
In particular, for sufficiently small $\ep>0$, there exists only one bounded positive solution, which in addition is periodic, hyperbolic attractive, and attracts all other positive solutions as time increases. That is, \eqref{eq:yeps} exhibits uniform stability.
\end{itemize}
\end{teor}
The fundamental ideas underlying the proof of Theorem \ref{teor:relaxation_oscillations} are twofold.
First, equation \eqref{eq:yeps} can be interpreted as the following slow-fast system, with derivatives expressed in terms of the fast time variable $t$:
\begin{equation}
\label{eq:fast_slow_system}
\left\{
\begin{array}{lll}
\theta'&=&\ep\,,\\[1ex]
z'&=&y^\ep(\theta)+\bar g(z)\,,
\end{array}
\right.
\end{equation}
where $\theta$ stands for the slow time variable and $y^\ep(\theta):=\alpha+(\beta+\ep^r)\sin\theta=y_\ep(\theta/\ep)$.
In this formulation, the application of Tikhonov's Theorem (see, e.g., Theorem 11.1 of Ref.~\onlinecite{khalil2002}) guarantees the tracking of the branches of hyperbolic critical points of the autonomous problem \eqref{eq:autonomousbonifacio}, provided that our solutions remain away from the autonomous bifurcation points. The second idea involves the construction a majorant equation and a minorant equation for \eqref{eq:yeps} in the neighborhood of each autonomous bifurcation point. The proof is long and technical, based on some ideas of Section 2.2 of Ref.~\onlinecite{berglund2006}. It is detailed in Appendix \ref{app}.

\begin{figure}
     \centering
         \includegraphics[width=\columnwidth]{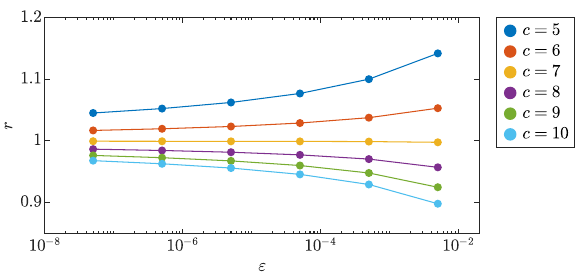}
         \caption{Numerical estimate (up to four decimal places) of the value of $r(\ep,c)$ at which the dynamical change of \eqref{eq:yeps} from relaxation solutions to bistability occurs, obtained using a bisection method following integration with MATLAB 2024b’s \texttt{ode45} solver over a period, starting from an initial condition very close to the upper stable branch at $t=0$.}
        \label{fig:simulationr}
\end{figure}

Fig.~\ref{fig:relaxation_simulation} depicts the scenarios of bistability and uniform stability appearing for different values of $\ep$ and $r$.
Fig.~\ref{fig:simulationr} presents a numerical approximation of the value of $r(\ep,c)$ such that the equation has relaxation solutions for $r>r(\ep,c)$ and exhibits bistability for $r<r(\ep,c)$, for various values of $c$ and small $\ep$. The convergence stated in Theorem \ref{teor:relaxation_oscillations} can be observed---albeit very slowly---as $\ep$ decreases, with the values tending toward $r = 1$.

\begin{figure}
     \centering
         \includegraphics[width=\columnwidth]{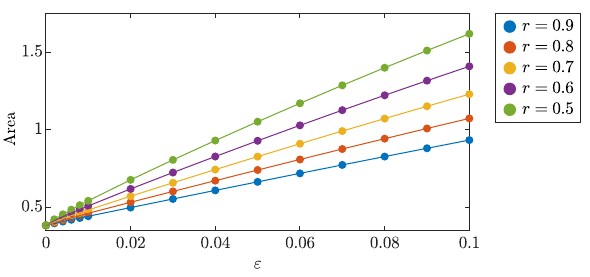}
         \caption{Numerical approximation of the area enclosed by the relaxation solution loop in \eqref{eq:yeps} for $c=5$ and various values of $\ep$ and $r$. MATLAB 2024b’s \texttt{polyarea} algorithm has been used.
         For all values of $r$, the area converges to that enclosed by the curve $\Gamma$ as $\ep\to0^+$.
         Linear fits of these curves performed using a power law of the form $\text{Area}=\text{Area}_\Gamma+C\,\ep^s$ exhibit $r$-dependent values of $C$ and $s$, where $\text{Area}_\Gamma$ stand for the area enclosed by $\Gamma$---approximated via \texttt{polyarea} from the autonomous bifurcation diagram.
}
        \label{fig:simulationareas}
\end{figure}

In physics, hysteresis refers to processes in which the state of a system depends on its past history---represented in this context by the input $y(t)$.
The extent of hysteresis is often characterized by the width or area of the hysteresis loop obtained when an appropriate system variable is plotted as a function of a switching parameter.\cite{hohl1995,jung1990}
In Figure \ref{fig:simulationareas}, the area enclosed by the hysteresis curve is numerically approximated for various values of $c$ and $\ep$.
\section*{Acknowledgements}
All the authors were supported by Ministerio de Ciencia, Innovaci\'{o}n y Universidades (Spain)
under project PID2021-125446NB-I00 and by Universidad de Valladolid under project PIP-TCESC-2020.
J.~Due\~{n}as was also supported by Ministerio de Universidades (Spain) under programme FPU20/01627.
\section*{Author declarations}
\vspace{-2ex}\noindent\textbf{Conflict of Interest}\par
The authors have no conflicts to disclose.\par\vspace{2ex}
\noindent\textbf{Author contributions}\par
\noindent\textbf{J.~Due\~{n}as:} Conceptualization, writing -- review \& editing.
\textbf{C.~N\'{u}\~{n}ez:} Conceptualization, writing -- review \& editing.
\textbf{R.~Obaya:} Conceptualization, writing -- review \& editing.
\section*{Data availability}
Data sharing is not applicable to this article as no new data were created or analyzed in this study.
\appendix
\section{Proof of Theorem \ref{teor:relaxation_oscillations}}\label{app}

Recall that the value of $r$ is fixed and $\ep$ varies. The core of the proof consists in analyzing the Poincar\'{e} map
\[
T_\ep\colon\R\to\R\,,\quad \bar x\mapsto x_\ep(2\pi/\ep,0,\bar x)\,,
\]
where $x_\ep(t,s,\bar x)$ stands for the solution of \eqref{eq:yeps} satisfying $x_\ep(s,s,\bar x)=\bar x$.
We assume from the beginning that $\ep<1$ satisfies $\lb_1(c)-\ep^r\ge 0$, and fix any $\rho>x_2(c)$ with $\bar g(x)\le-\lb_2(c)-1$ for $x\ge \rho$, so the graph on any solution of \eqref{eq:yeps} eventually enters in the set $\R\times [0,\rho]$, which is forward invariant and contains the graphs of all the bounded solutions.
\smallskip

(i) Let us prove the claim for $r>1$.
Our goal is to check three properties for sufficiently small $\ep>0$: (A) that the restriction of $T_\ep$ to the interval $[x_2(c),\rho]$ is well-defined and has a globally attractive fixed point within $[x_2(c),\rho]$, corresponding to a unique attractive hyperbolic (periodic and positive) solution of \eqref{eq:yeps}; (B) the same statement, but applied to the interval $[0,\,x_1(c)]$; and (C) that the restriction of the inverse map $T^{-1}_\ep$ of $T_\ep$ to the interval $[x_1(c),\,x_2(c)]$ is well-defined and has a unique globally attractive fixed point in $[x_1(c),\,x_2(c)]$, corresponding to a repulsive hyperbolic solution of \eqref{eq:yeps}. In this situation, which clearly precludes the occurrence of relaxation solutions, bistability can be proved by repeating part of the arguments of the proof of Theorem \ref{teor:bistability_comparison_autonomous}.

We prove statement (A) in three steps.

\textbf{Step 1.} \textit{Tracking the uniformly stable branch.}
We consider the slow-fast system \eqref{eq:fast_slow_system}. Let $(\bar\theta+\ep\,t,z^\ep(t,\bar z,\bar\theta))$ be the solution of \eqref{eq:fast_slow_system} with value $(\bar\theta,\bar z)$ at $t=0$.
Then, $z^\ep(t,\bar\theta,\bar x)=x_\ep(t+\bar\theta/\ep,\bar\theta/\ep,\bar x)$.
For $\ep=0$, we get $z'=y^0(\theta)+\bar g(z)$ with $y^0(\theta):=\alpha+\beta\sin\theta$, which takes values between $\lb_1(c)$ and $\lb_2(c)$. Let $u(\lb)$ be the upper stable equilibrium of \eqref{eq:autonomousbonifacio}$_\lb$ for $\lb>\lb_1(c)$: $\lambda\mapsto u(\lb)$ is the upper red solid curve in Fig.~\ref{fig:bif_diagram_aut_bonifacio}. Then, $y^0(\theta)+\bar g(u(y^0(\theta)))=0$ for all $\theta\in[-\pi/2,\,3\pi/2]$. We fix $\lb_0\in(\lb_1(c),\,\alpha)$, define $\mD_0:=\{\theta\in[-\pi/2,\,3\pi/2]\mid\,y^0(\theta)\geq\lb_0\}$, and check that $-\pi/2<\inf\mD_0<0<\pi<\theta_0:=\sup\mD_0<3\pi/2$.
It is not hard to check the existence of $k\geq 1$ and $\beta>0$ such that $|z^0(t,\bar\theta,\bar z)-z^0(t,\bar\theta,u(y^0(\bar\theta)))|=|z^0(t,\bar\theta,\bar z)-u(y^0(\bar\theta))|\leq k\,e^{-\beta\,t}|\bar z-u(y^0(\bar\theta))|$ for all $\bar\theta\in\mD_0$, $\bar z\in[x_2(c),\,\rho]$ and $t\geq0$.
So, Tikhonov's Theorem (see, e.g., Theorem 11.1 of Ref.~\onlinecite{khalil2002}) ensures that $x_\ep(\theta_0/\ep,0,\bar x)-u(y^0(\theta_0))=O(\ep)$ uniformly for $\bar x\in[x_2(c),\,\rho]$ if $\ep>0$ is sufficiently small.

\textbf{Step 2.} \textit{Bypassing the autonomous bifurcation point.}
We will prove that, if $\ep>0$ is sufficiently small, then $\tilde x_\ep(t):=x_\ep(t,0,x_2(c))$ satisfies $\tilde x_\ep(t)>x_2(c)$ for all $t\in(0,\,2\pi/\ep]$.

Let us call $t_\ep^-<t_\ep^+$ the two points of $[0,2\pi/\ep]$ with $y_\ep(t_\ep^\pm)=\lb_1(c)$. It is easy to check that
$\theta_0/\ep<t_\ep^-<3\pi/(2\ep)<t_\ep^+<2\pi/\ep$ if $\ep>0$ is sufficiently small. Since $\tilde x'_\ep(0)>0$ and $y_\ep(t)+\bar g(x_2(c))=y_\ep(t)-\lb_1(c)$, an easy contradiction argument shows that: $\tilde x_\ep(t)>x_2(c)$ for all $t\in(0,\,t_\ep^-)$; and, if $\tilde x_\ep(t_\ep^+)>x_2(c)$ (as we will prove below), then  $\tilde x_\ep(t)>x_2(c)$ for all $t\in[t_\ep^+,\,2\pi/\ep]$\,.
Hence, Step 2 will be completed once it has been proved that, if $\ep>0$ is sufficiently small, then $\tilde x_\ep(t)>x_2(c)$ for $t\in[t_\ep^-,\,t_\ep^+]$. In turn, since $u(y^0(\theta_0))=u(\lb_0)$, Step 1 ensures that $\tilde x_\ep(\theta_0/\ep)>(u(\lb_0)+x_2(c))/2$ if $\ep>0$ is sufficiently small (what we assume from now on), so it suffices to prove that  $\bar x_\ep(t):=x_\ep(t,\theta_0/\ep,(u(y^0(\theta_0))+x_2(c))/2)>x_2(c)$ for $t\in[t_\ep^-,\,t_\ep^+]$. Observe also that there is no restriction in assuming that $(u(y^0(\theta_0))+x_2(c))/2<\rho$.

Note that $y_\ep(t)\ge y_\ep(t_\ep^-)+y_\ep'(t_\ep^-)(t-t_\ep^-)$ on $[\pi/\ep,2\pi/\ep]$, where $y_\ep$ is strictly convex. It is not hard to check that
\[
\cos(\ep\,t_\ep^-)=-\sqrt{\frac{\frac{\ep^r}{\beta}\left(\frac{\ep^r}{\beta}+2\right)}{\left(\frac{\ep^r}{\beta}+1\right)^2}}
=-\sqrt{2/\beta}\,\ep^{r/2} + O(\ep^{3r/2})
\]
as $\ep\to 0^+$. So,
\[
\begin{split}
y_\ep'(t_\ep^-)\,(t-t_\ep^-)&=\ep\,(\beta+\ep^r)\cos(\ep\,t_\ep^-)\,(t-t_\ep^-)\\
&\ge -k_1\;\ep^{r/2+1}(t-t_\ep^-)-k_2\,\ep^{3r/2+1}
\end{split}
\]
for $t\in[\pi/\ep,\,2\pi/\ep]$ if $\ep$ is small, with $k_1:=\sqrt{2\beta}$ and $k_2>0$.

Note also that, since $\bar g'(x_2(c))=0$, Taylor's Theorem shows that, if $x\in[x_2(c),\,\rho]$, then $\bar g(x)\ge \bar g(x_2(c))-k_3\,(x-x_2(c))^2$ for $-k_3:=\inf_{x\in[x_2(c),\,\rho]}\bar g''(x)/2<0$.
So, $\bar x_\ep(t)$ satisfies
\[
\begin{split}
 x'&\ge y_\ep(t_\ep^-)+y_\ep'(t_\ep^-)(t-t_\ep^-)+\bar g(x_2(c))-k_3\,(x-x_2(c))^2\\
 &=y_\ep'(t_\ep^-)(t-t_\ep^-)-k_3\,(x-x_2(c))^2\\
 &\ge-k_1\,\ep^{r/2+1}(t-t_\ep^-)-k_2\,\ep^{3r/2+1}-k_3\,(x-x_2(c))^2
\end{split}
\]
as long as $t\in[\pi/\ep,\,2\pi/\ep]$ and $\bar x_\ep(t)\ge x_2(c)$.
We make the double change of variables
\begin{equation}
\label{eq:cambio}
\begin{split}
s&:=\ep^{r/6+1/3}(t-t_\ep^-)\,,\\
w&:=\ep^{-r/6-1/3}(x-x_2(c))\,,
\end{split}
\end{equation}
chosen to transform the previous differential inequality in
\begin{equation}
\label{eq:simpleRiccatiequation}
\dot w\geq-\beta_\ep-k_1\,s-k_3\,w^2
\end{equation}
for $\beta_\ep:=k_2\,\ep^{7r/6+1/3}$, where $\dot w:=dw/ds$.
This change transforms $\bar x_\ep$ into a map $\bar w_\ep$ that takes the value
\[
 w_0^\ep:=\ep^{-r/6-1/3}(u(y^0(\theta_0))-x_2(c))/2
\]
at the time $s_0^\ep:=\ep^{r/6-2/3}(\theta_0-\ep\,t_\ep^-)<0$, which in addition solves the inequation \eqref{eq:simpleRiccatiequation} as long as $s\in[\ep^{r/6-2/3}(\pi-\ep\,t_\ep^-),$ $\,\ep^{r/6-2/3}(2\pi-\ep\,t_\ep^-)]$ and $\bar w_\ep(s)\ge 0$.

The solution $w(t)$ of $\dot w=-k_1\,s-k_3\,w^2$ with $w(0)=0$ satisfies $w'(0)=0$ and $w''(0)=-1$, so that it is strictly negative for small $s\neq0$ and strictly concave on an interval containing $0$. By continuous dependence, there exists $w_0>0$ such that the solution $\tilde w(t)$ of $\dot w=-k_1\,s-k_3\,w^2$ with $\tilde w(0)=w_0$ is strictly positive on an interval $(s^*_1,s^*_2)\ni 0$ and vanishes at its endpoints. We take a $\ep>0$ small enough to get $w_0^\ep\ge\max_{s\in[s^*_1,\,s^*_2]}\tilde w(s)$. It is obvious that, if $s_0^\ep\in[s_1^*,\,0)$, then the solution of $\dot w=-k_1\,s-k_3\,w^2$ taking the value $w_0^\ep$ at $s_0^\ep$ is strictly positively bounded from below on $[s^*_1,s^*_2]$. The same happens if $s_0^\ep<s_1^*$, since the same solution is above that of $\dot w=-k_1\,s-k_3\,w^2$ with value $0$ at $s_0^\ep$. Taking a possible smaller $\ep>0$, we conclude that the solution of $\dot w=-\beta_\ep-k_1\,s-k_3\,w^2$ with the same initial data is strictly positive on $[s^*_1,s^*_2]$. Retracing the steps, we conclude that $\bar x_\ep(t)$ is strictly above $x_2(c)$ at least on $[t_\ep^-,\,t_\ep^-+\ep^{-r/6-1/3}s_2^*]$. The proof will be completed once it has been checked that $t_\ep^-+\ep^{-r/6-1/3}s_2^*\ge t_\ep^+$ if $\ep>0$ is sufficiently small.

It is easy to check that $\cos((\ep\,t_\ep^+-\ep\,t_\ep^-)/2)=\cos(3\pi/2-\ep\,t_\ep^-)=-\sin(\ep\,t_\ep^-)=1/(1+\ep^r/\beta)$. Since $\arccos(1/(1+x))=\sqrt{2\,x}+O(x^{3/2})$ as $x\to0^+$ (as $\lim_{x\to 0^+}(\arccos(1/(1+x))-\sqrt{2\,x})/x^{3/2}= -5\sqrt{2}/12$ ensures), we get
\begin{equation}\label{eq:landau_sep-tep}
\ep\,(t_\ep^+-t_\ep^-)=2\sqrt{2/\beta}\,\ep^{r/2}+O\big(\ep^{3r/2}\big)\,.
\end{equation}
Since $r>1$, $\lim_{\ep\to 0^+}\ep^{r/2-1}/\ep^{-r/6-1/3}=\lim_{\ep\to 0^+}\ep^{2(r-1)/3}=0$, so $t_\ep^+-t_\ep^-<\ep^{-r/6-1/3}s_2^*$ if $\ep>0$ is sufficiently small. (For next purposes, note that this last is the only point in the proof in which the condition $r>1$ plays a role.) Step 2 is complete.

\textbf{Step 3. }\textit{Unique hyperbolic fixed point of $T_\ep$.}
The choice of $\rho$ and the result of the previous step yields $x_2(c)<\tilde x_\ep(t)\le x_\ep(t,0,\bar x)\leq\rho$ for all $t\in(0,\,2\pi/\ep]$ if $\bar x\in[x_2(c),\,\rho]$. Hence, the restriction of $T_\ep$ to $[x_2(c),\rho]$ is well defined. In addition, since $\bar g'$ is strictly negative on $(x_2(c),\rho]$,
\[
 T_\ep'(\bar x)-1=\!\int_0^{2\pi/\ep}\!\bar g'(x_\ep(t,0,\bar x))\,\exp\left(\!\int_0^t\bar g'(x_\ep(s,0,\bar x))\,ds\!\right)dt
\]
is also strictly negative for all $\bar x\in[x_2(c),\,\rho]$. Since $T'_\ep(\bar x)=\exp\big(\int_0^{2\pi/\ep}\bar g'(x_\ep(s,0,\bar x))\,ds\big)$ is continuous and positive, we conclude that $|T_\ep'(\bar x)|\le k<1$ for all $\bar x\in (x_2(c),\rho]$, and Banach's Theorem ensures that the restriction is contractive: there exists a unique fixed point of $T_\ep$ on $[\bar x_2(c),\rho]$ and it is globally hyperbolic attractive---therefore corresponding to the unique (attractive hyperbolic) periodic solution of \eqref{eq:yeps} on $[x_2(c),\,\infty)$.\hspace{-1cm}~
\smallskip

This completes the proof of (A). The arguments used to prove (B) and (C) are analogous to those already detailed. In particular, those used to  prove (B) are entirely symmetrical. To prove (C), we work with the time-reversed equation and apply twice the arguments employed in the vicinity of the autonomous bifurcation point presented in Step 2, as one of them appears at each end of the branch of repulsive equilibria tracked by the hyperbolic solution. These arguments complete the proof of (i).
\smallskip

(ii) Let us now assume $r<1$. We take $t_\ep^-$ as in the proof of (i), and determine $\bar t_\ep^-\in(0,\,\pi/(2\ep))$ by $y_\ep(\bar t_\ep^-)=\lb_2(c)$. Our first claim is: if $\ep>0$ is sufficiently small, then $\tilde x_\ep(t):=x_\ep(t,0,\bar x)$ with $\bar x\in[0,\,\rho]$ satisfies  $\tilde x_\ep(t)>x_2(c)$ for all $t\in(\pi/(2\ep),\,t_\ep^-)$ and $\tilde x_\ep(t)<x_1(c)$ for all $t\in(3\pi/(2\ep),\,2\pi/\ep+\bar t_\ep^-]$. And, if $\bar x\in[0,\,x_1(c)]$, then $\tilde x_\ep(t)<x_1(c)$ for all $t\in(0,\,\bar t_\ep^-]$. Assume that it is true. Then, $T_\ep\colon[0,\rho]\to[0,\,x_1(c)]$ is well defined. In addition, $T_\ep'(x)<1$ for all $x\in[0,\,x_1(c)]$, as we can prove by making the change of variables $l=\ep\,t$ on the expression of $T_\ep'(x)-1$ (see Step 3 of (i)) and observing that the intervals on which the (bounded) integrand is non-positive are contained in $[\ep\,\bar t_\ep^-,\,\pi/2]\cup [\ep\,t_\ep^-,\,3\pi/2]$ and that the length of this union has limit $0$ as $\ep\to 0^+$. So, the unique fixed point of $T_\ep$ on $[0,\,x_1(c)]$ determines the unique periodic solution of \eqref{eq:yeps}, which is hyperbolic attractive and with global domain of attraction. This proves the uniform stability.

Assume for the moment being that $\tilde x_\ep(t_0)\ge x_2(c)$ for a point $t_0\in (\bar t_\ep^-,\,\pi/\ep)$ (resp.~$\tilde x_\ep(t_0)\le x_1(c)$ for $t_0\in(t_\ep^-,\,2\pi/\ep)$). Reasoning as in Step 2 of (i) we check that $\tilde x_\ep(t)>x_2(c)$ for all $t\in(t_0,\,t_\ep^-)$ (resp.~$\tilde x_\ep(t_\ep^+)\le x_1(c)$ for all $t\in(t_0,\,2\pi/\ep+\bar t_\ep^-)$). The same argument shows that then $\tilde x_\ep(t)<x_1(c)$ for all $t\in(0,\,\bar t_\ep^-]$ if $\bar x\in[0,\,x_1(c)]$. With this and basic comparison arguments in mind, it is clear that our claim will be proved once seen that $x_\ep(t,t_\ep^-,\,\rho)$ reaches $x_1(c)$ before $3\pi/(2\ep)$ and that $x_\ep(t,\bar t_\ep^-,0)$ reaches $x_2(c)$ before $\pi/(2\ep)$. Let us focus on the first assertion: the second one is proved similarly.

First, we fix a point $x_*\in(\sqrt{3},\,x_2(c))$. Let $s_*$ be any time in $(t_\ep^-,\,3\pi/(2\ep))$. It easy to check the existence of a constant $k>0$ such that the solution $x_\ep(t,s_*,\bar x)$ with $\bar x\in(x_1(c),x_*]$ satisfies $x'\le -k$ as long as $t\in(s,\,3\pi/2\ep)$ and $x_\ep(t,s_*,\bar x)\ge x_1(c)$. So, $x_\ep(t,s_*,\bar x)$ reaches $x_1(c)$ at a time not greater than $t_1:=(x_*-x_1(c))/k$, and we can assume that $\ep>0$ is small enough as to get $t_1\le(3\pi/(2\ep)-t_\ep^-)/2$ since $r<1$. So, our goal is reduced to check that $x_\ep(t,t_\ep^-,\,\rho)$ reaches $x_*$ before $3\pi/(2\ep)-t_1$.

Taylor's Theorem provides $k_5>0$ with $\bar g(x)\le-\lb_1(c)-k_5\,(x-x_2(c))^2$ for all $x\in[x_*,\rho]$. On the other hand, the convexity of $y_\ep$ on $[t_\ep^-,\,3\pi/(2\ep)]$ ensures that its graph over $[t_\ep,3\pi/(2\ep)]$ is below the line joining the points $(t_\ep,\lb_1(c))$ and $(3\pi/(2\ep),\lb_1(c)-\ep^r)$; i.e.,
$y_\ep(t)\le\lb_1(c)-\ep^r\,(3\pi/(2\ep)-t_\ep^-)^{-1}(t-t_\ep^-)$ for $t\in[t_\ep^-,\,3\pi/(2\ep)]$, which combined with $(3\pi/(2\ep)-t_\ep^-)=(t_\ep^+-t_\ep^-)/2$ and \eqref{eq:landau_sep-tep} ensures the existence of $k_4>0$ such that $y_\ep(t)\le\lb_1(c)-k_4\ep^{r/2+1}(t-t_\ep^-)$ for $t\in[t_\ep^-,\,3\pi/(2\ep)]$ if $\ep>0$ is sufficiently small. Therefore, we can assume that $\bar x_\ep(t):=x_\ep(t,t_\ep^-,\rho)$ satisfies
\[
 x'\le -k_4\,\ep^{r/2+1}\,(t-t_\ep^-)-k_5\,(x-x_2(c))^2
\]
as long as $t\in[t_\ep^-, 3\pi/(2\ep)]$ and $x\ge x_*$. The changes of variables \eqref{eq:cambio} take this inequality to
\[
 \dot w\le -k_4\,s-k_5\,w^2
\]
and transform $\bar x_\ep$ into a forward bounded map $\bar w_\ep$ with $\bar w_\ep(0)=\ep^{-r/6-1/3}(\rho-x_2(c))$ that satisfies the previous inequation as long as $s\in[0,\,\ep^{r/6-2/3}(3\pi/2-\ep\,t_\ep^-)]$ and $\bar w_\ep(s)\ge \ep^{-r/6-1/3}(x_*-x_2(c))$. It is not hard to check that $\bar x'_\ep$ does not vanish on $[t_\ep^-,\,3\pi/(2\ep)]$, and hence $\dot{\bar{w}}_\ep$ does not vanish on $[0,\,\ep^{r/6-2/3}(3\pi/2-\ep\,t_\ep^-)]$. The new change of variable $w=\dot z/(k_5\,z)$ takes $\dot w=-k_4\,s-k_5\,w^2$ into $\dot{\dot z}+k_4\,k_5\,s\,z=0$. Applying Sturm's Separation and Comparison Theorems (see, e.g., Corollary XI.3.1 of Ref.~\onlinecite{hartman2002}) to the majorant $\dot{\dot z}+k_4\,k_5\,z=0$ on $[1,\infty)$ and observing that the map $\sin(\sqrt{k_4k_5}\,s)$ solves this last equation, we conclude that any solution of $\dot w=-k_4\,s-k_5\,w^2$ vanishes at least once at $[t_2/2,\,t_2]$ for $t_2:=2\pi/\sqrt{k_4k_5}$. An easy contradiction argument shows the existence of $s_\ep\in [0,\,2\pi]$ such that $\bar w_\ep(s_\ep)=\ep^{-r/6-1/3}(x_*-x_2(c))$, which means that $\bar x_\ep(t_\ep^-+\ep^{-r/6-1/3} s_\ep)=x_*$. Altogether, it suffices to check that $t_\ep^-+\ep^{-r/6-1/3} t_2\le 3\pi/(2\ep)-t_1$ if $\ep>0$ is sufficiently small, and this is made using again $(3\pi/(2\ep)-t_\ep^-)=(t_\ep^+-t_\ep^-)/2$ and \eqref{eq:landau_sep-tep}, and applying the fundamental condition $r<1$. The proof of our uniform stability is hence complete.

Let $x^*_\ep(t)$ be the $2\pi/\ep$-periodic solution of \eqref{eq:yeps} for $\ep>0$ sufficiently small. It remains to prove that the $\ep$-family $x^*_\ep(t)$ of \eqref{eq:yeps} is a relaxation oscillation. By reviewing the proof, we observe that $x^*_\ep(t)\le x_1(c)$ for $t\in[0,\,\bar t_\ep^-]\cap [t_\ep^+,2\pi/\ep]$ and $x^*_\ep(t)\ge x_1(c)$ for $t\in[\bar t_\ep^+,\,t_\ep^-]$: see Fig.~\ref{fig:relaxation_simulation}.
Recall also that $\bar t_\ep^\pm$ and $t_\ep^\pm$ are the points at which $y_\ep$ crosses the autonomous bifurcation values $\lb_2(c)$ and $\lb_1(c)$. The fact that $\{(y_\ep(t),\,x^*_\ep(t))\mid\,t\in\R\}$ converges in Hausdorff distance to the curve $\Gamma$ as $\ep\to0^+$ is proved by combining these properties with a careful application of Tikhonov's Theorem. We do not include the details.

\bibliography{CCApplicationsBib}

\end{document}